\documentclass[5p]{elsarticle}

\usepackage[intlimits]{amsmath}
\usepackage{amsthm,amssymb}
\usepackage{color}

\usepackage[english]{babel}
\usepackage[]{algorithm2e}
\SetKwProg{Init}{Initialization:}{}{}

\usepackage{graphics} 
\usepackage{epsfig} 
\theoremstyle{plain}
\newtheorem{assumption}{Assumption}

\newtheorem{theorem}{Theorem}
\newtheorem{proposition}{Proposition}

\theoremstyle{remark}
\newtheorem{remark}{Remark}
\biboptions{sort&compress}
\usepackage{url}

\journal{Systems \& Control Letters}

\begin{document}
	
	\begin{frontmatter}
		
		\title{On Computation of Approximate Solutions to Large-Scale Backstepping Kernel 
		Equations via Continuum Approximation\tnoteref{t1}}
		\tnotetext[t1]{Funded by the European Union (ERC, C-NORA, 101088147). Views and 
		opinions 
		expressed are however those of the authors only and do not necessarily reflect those of the 
		European Union or the European Research Council Executive Agency. Neither the European 
		Union nor the granting authority can be held responsible for them.}
		
		\author[1]{Jukka-Pekka Humaloja}
		\author[1]{Nikolaos Bekiaris-Liberis}
		
		\affiliation[1]{organization={Department of Electrical and Computer Engineering, Technical 
		University of Crete},
			addressline={University Campus, Akrotiri}, 
			city={Chania},
			postcode={73100}, 
			country={Greece}}
		
		\begin{abstract}
We provide two methods for computation of continuum backstepping kernels that arise in 
control 
of continua (ensembles) of linear hyperbolic PDEs and which can approximate backstepping 
kernels arising in control of a large-scale, PDE system counterpart (with computational 
complexity 
that does not grow with the number of state components of the large-scale system). In the 
first method, we provide explicit formulae for the solution to the continuum 
kernels PDEs, employing a (triple) power series representation of the continuum kernel and 
establishing its convergence properties. In this case, we also provide means for reducing 
computational complexity by properly truncating the power series (in the powers of the 
ensemble variable). In the second method, we identify a class of systems for which the solution to 
the continuum (and hence, also an approximate solution to the respective large-scale) kernel 
equations can be constructed in closed form. We also present numerical examples to illustrate 
computational efficiency/accuracy of the approaches, as well as to validate the stabilization 
properties of the approximate control kernels, constructed based on the continuum.  
		\end{abstract}
		
		
	\end{frontmatter}

\section{Introduction}

Exponentially stabilizing feedback control laws  based on the backstepping method 
\cite{KrsSmyBook} have been developed for several classes of hyperbolic PDEs in recent years, 
see
\cite{AllKrs23arxiv, DiMVaz13, YuHKrs21,AurBre22, RedAur21, CorVaz13,HuLDiM16,WanKrs22, 
DiaDia17}. The application of the method requires 
computing the backstepping control (and observer) gains, which involves solving the respective 
kernel PDEs. In certain instances, as in the case of $2\times2$ linear, hyperbolic systems with 
constant coefficients \cite{VazKrs14}, it is possible 
to derive the solution to the kernel PDEs in closed-form; whereas, in general, the solution to the 
kernel PDEs has to be computed (i.e., approximated) numerically \cite{AnfAamBook, 
VazCheCDC23,AurMor19,QiJZha24,BhaShi24}. However, for large-scale hyperbolic systems,  
such numerical schemes 
for solving the kernel PDEs exhibit computational complexity that grows with the number of state 
components of the PDEs \cite{HumBek24arxiv, HumBekCDC24}. Motivated by this, in the 
present paper, 
we provide tractable computational tools for constructing exponentially stabilizing feedback laws 
(by solving the respective kernel PDEs) for large-scale \cite{DiMVaz13} and continua 
\cite{AllKrs23arxiv} of hyperbolic PDEs. We achieve this in two ways; by adapting the power series 
method \cite{VazCheCDC23} to solve continuum kernel PDEs and by providing closed-form 
(exact) solutions for a class of such kernel PDEs. As we have shown in 
 \cite{HumBek24arxiv, HumBekCDC24}, the  continuum control kernel provides stabilizing 
 feedback control gains 
for the respective large-scale hyperbolic PDE system as well.

For large-scale, 1-D, linear hyperbolic systems, one of the most efficient approaches to 
compute backstepping control kernels is to utilize a closed-form, explicit solution. However, an 
exact, 
closed-form solution is available only for specific classes of hyperbolic systems, such as, for 
$2\times2$ systems with constant coefficients \cite{VazKrs14}. An alternative in providing explicit, 
backstepping control kernels, constitutes in computing approximate (yet explicit) solutions to the 
respective kernels PDEs via, for example, late lumping \cite{AurMor19}, neural operators 
\cite{BhaShi24}, and power series \cite{VazCheCDC23, AscAst18, VazZha23, LinVaz24arxiv} 
-based approaches. In particular, the power series method provides a simple and flexible 
approach for computation of approximate backstepping kernels, thus making it a suitable 
candidate for computation of the solutions to continua and large-scale kernels PDEs.

As the main contribution of the paper, we present a triple power series-based approach to 
compute the 
solution to continuum kernel equations on a prismatic 3-D domain, thus, providing a systematic, 
computational tool for constructing stabilizing feedback laws for continua of hyperbolic PDEs. 
In turn, utilizing our recent results on continuum approximation of large-scale systems of 
hyperbolic PDEs   \cite{HumBek24arxiv, HumBekCDC24}, the stabilizing feedback laws 
obtained for the continuum  
can be directly employed for stabilizing the corresponding large-scale system.
We provide theoretical guarantees that the power series approximation is convergent, provided 
that the parameters of the continuum kernel PDEs are analytic, and present an algorithm to 
compute the solution. Moreover, we propose an order-reduction method for the power series 
approach, which can potentially reduce the computational complexity of the approach from 
$\mathcal{O}(N^3)$ to $\mathcal{O}(N^2)$, where $N$ is the order of the power series 
approximation. In either case, we show that the power series method to compute continuum 
kernels results in providing stabilizing control kernels for the respective large-scale system with 
computational complexity that does not grow with the number of PDE state components (in 
contrast to the case of employing the power series approach to solve the large-scale kernel 
equations). We then apply the proposed approach to a stabilization problem of a large-scale 
hyperbolic system of PDEs, where an exponentially stabilizing state-feedback gain can be 
computed efficiently by combining the continuum approximation approach  
\cite{HumBek24arxiv, HumBekCDC24} 
and the proposed reduced-order, triple power series method. We also conduct numerical 
experiments to illustrate the computational effectiveness/accuracy of the (triple) power series 
method.

 We then identify special cases, in which the solution to the continuum kernel equations can be 
 expressed in closed form. This may be viewed as reminiscent of \cite{VazKrs14}, where explicit 
 kernels are computed for stabilization of $2\times 2$ linear hyperbolic PDE systems with 
 constant 
 parameters. This special class of continuum kernel equations is characterized by specific 
 conditions on the respective continuum parameters, under which we provide a closed-form 
 solution in a constructive manner. Because the conditions imposed concern the continuum 
 parameters, there is some flexibility degree for their satisfaction in the case in which the 
 continuum parameters involved are obtained as continuum approximation (which may not be 
 unique) of the respective sequences of parameters of the large-scale system counterpart (see 
 \cite{HumBek24arxiv, HumBekCDC24}). We also present a numerical example for 
 which these conditions on the continuum kernel 
 equations are satisfied, while a closed-form solution cannot be derived for the respective 
 large-scale system. Thus, such closed-form solution provides explicit stabilizing control gains for 
 the continuum, as well as for any large-scale system of hyperbolic PDEs that can be
approximated by such a continuum, even when a closed-form solution for the original, large-scale 
kernel equations may not be available (see also  \cite{HumBek24arxiv, HumBekCDC24}).

The paper is organized as  follows. In Section~\ref{sec:appr}, we present the large-scale kernel 
equations (associated to a large-scale system of hyperbolic PDEs) and their continuum 
approximation. In 
Section~\ref{sec:pow}, we present the power series approach to solving the continuum kernel 
PDEs, including an explicit algorithm and a potential order-reduction method.  In 
Section~\ref{sec:exp}, 
we identify sufficient conditions for the parameters of the continuum kernel equations, under 
which the solution can be expressed in closed form. In Section~\ref{sec:ex1}, 
we demonstrate the accuracy and convergence rate of the power series approximation in a 
numerical example. In Section~\ref{sec:ex2}, we 
apply the power series method to the stabilization problem of a large-scale system of hyperbolic 
PDEs, where a closed-form solution to the kernel equations (large-scale or continuum) is not 
available. Finally, Section~\ref{sec:conc} contains concluding remarks.

\section{Large-Scale Kernel Equations and Their Continuum Approximation} \label{sec:appr}

In this section, we present the considered class of large-scale hyperbolic kernel equations, 
and how they can be approximated by respective continua of kernel equations by virtue of our 
earlier works \cite{HumBek24arxiv, HumBekCDC24}. This provides the basis for the development 
of computational tools to approximate the solution of the large-scale kernel equations by solving 
the respective continua of kernel equations instead.

For $i=1,\ldots,n$, where $n$ is large, consider kernel equations of the form 
\cite{DiMVaz13, 
HumBek24arxiv}
\begin{subequations}
	\label{eq:kn}%
	\begin{align}
		\mu(x)\partial_xk^i(x,\xi) - \lambda_{i}(\xi)\partial_{\xi}k^i(x,\xi) & = \nonumber  \\
		\lambda_i'(\xi)k^i(x,\xi)+\frac{1}{n}\sum_{j=1}^n \sigma_{j,i}(\xi)k^j(x,\xi) +
		\theta_i(\xi)k^{n+1}(x,\xi),& \\
		\mu(x)\partial_xk^{n+1}(x,\xi) + \mu(\xi)\partial_{\xi}k^{n+1}(x,\xi) & = \nonumber \\
		-\mu'(\xi)k^{n+1}(x,\xi) +  \frac{1}{n}\sum_{j=1}^nW_j(\xi)k^j(x,\xi), & \label{eq:knb} 
	\end{align}
\end{subequations}
where $\mu > 0 , \lambda_i > 0, \sigma_{i,j}, \theta_i, W_i$, and $q_i$ are given parameters 
and 
$(k^i)_{i=1}^{n+1}$ is the sought solution, i.e., the exact kernels. The equations \eqref{eq:kn} are 
defined on the triangular domain
\begin{equation}
	\mathcal{T} = \left\{ (x,\xi) \in [0,1]^2: 0 \leq \xi \leq x \leq 1\right\},
\end{equation}
and equipped with  boundary conditions
\begin{subequations}
	\label{eq:knbc}%
	\begin{align}
		k^i(x,x) & = -\frac{\theta_i(x)}{\lambda_i(x) + \mu(x)},  \\
		\mu(0)k^{n+1}(x,0) & =  \frac{1}{n}\sum_{j=1}^nq_j\lambda_j(0)k^j(x,0),
	\end{align}
\end{subequations}
for all $x \in [0,1]$. Our aim is to find approximate solutions to \eqref{eq:kn}, \eqref{eq:knbc} by 
solving the corresponding continuum kernel equations of the form \cite{AllKrs23arxiv, 
HumBek24arxiv}
\begin{subequations}
	\label{eq:kc}%
	\begin{align}
		\mu(x)\partial_xk(x,\xi,y) - \lambda(\xi,y)\partial_{\xi}k(x,\xi,y) -  \theta(\xi,y)\bar{k}(x,\xi) & = 
		\nonumber \\
	k(x,\xi,y) \partial_{\xi}\lambda(\xi,y)
		+\int\limits_0^1\sigma(\xi,\eta,y)k(x,\xi,\eta)d\eta, & \label{eq:kca} \\
		\mu(x) \partial_x\bar{k}(x,\xi) + \mu(\xi)\partial_\xi\bar{k}(x,\xi) & =  \nonumber \\
		-\mu'(\xi)\bar{k}(x,\xi) + \int\limits_0^1W(\xi,y)k(x,\xi,y)dy, &
	\end{align}
\end{subequations}
where $\mu >0, \lambda>0, \sigma, \theta, W$, and $q$ are given parameters and $(k, 
\bar{k})$ is 
the sought solution, i.e., the continuum kernels. The equations \eqref{eq:kc} are defined on a 
prismatic domain
\begin{equation}
	\mathcal{P} = \left\{(x,\xi,y) \in [0,1]^3: (x,\xi) \in \mathcal{T}\right\},
\end{equation}
and equipped with boundary conditions 
\begin{subequations}
	\label{eq:kcbc}%
	\begin{align}
		k(x,x,y) & = - \frac{\theta(x,y)}{\lambda(x,y) + \mu(x)}, \label{eq:kcbca} \\
		\mu(0)\bar{k}(x,0) & = \int\limits_0^1q(y)\lambda(0,y)k(x,0,y)dy, \label{eq:kcbcb}
	\end{align}
\end{subequations}
for all $x \in [0, 1]$ and for almost every $y\in [0,1]$.

The kernel equations \eqref{eq:kn}, \eqref{eq:knbc} appear in backstepping control of $n+1$ 
hyperbolic PDEs \cite{DiMVaz13} (see \eqref{eq:n+1}, \eqref{eq:nbcuy} in 
\ref{app:n+1}), where the solution $\left(k^i\right)_{i=1}^{n+1 }$ provides an 
exponentially stabilizing state feedback gain for the system (see \eqref{eq:Un}). More recently 
\cite{AllKrs23arxiv}, 
the backstepping control methodology has been extended to continua of hyperbolic PDEs (see 
\eqref{eq:inf}, \eqref{eq:cbcuy} in \ref{app:n+1}), where 
the exponentially stabilizing state feedback gain is obtained from the solution $(k,\bar{k})$ to the 
continuum kernel equations \eqref{eq:kc}, \eqref{eq:kcbc} (see \eqref{eq:Uc}). In our recent work  
\cite{HumBek24arxiv, HumBekCDC24}, we have shown that the continuum kernel equations 
\eqref{eq:kc}, 
\eqref{eq:kcbc} can be interpreted as an approximation of \eqref{eq:kn}, \eqref{eq:knbc} when 
$n$ is sufficiently large (see \cite[Lem. 4.2 \& Lem. 4.3]{HumBek24arxiv}). Thus, the 
exponentially 
stabilizing state feedback gain for the 
large-scale $n+1$ hyperbolic PDE can be approximated from the solution to  \eqref{eq:kc}, 
\eqref{eq:kcbc}, without compromising the exponential stability of the closed-loop system 
(see \cite[Thm 4.1]{HumBek24arxiv}). The 
relation between the parameters of \eqref{eq:kn}, \eqref{eq:knbc} and \eqref{eq:kc}, 
\eqref{eq:kcbc} was established in \cite{HumBek24arxiv, HumBekCDC24} as
\begin{subequations}
	\label{eq:afn1}%
	\begin{align}
		\lambda(x,i/n) & = \lambda_i(x), \\
		W(x, i/n) & = W_i(x), \\
		\theta(x,i/n) & = \theta_i(x), \\
		\sigma(x,i/n,j/n) & = \sigma_{i,j}(x), \label{eq:afns}\\	
		q(i/n) & = q_i,
	\end{align}%
\end{subequations}
for all $x \in [0,1]$ and $i,j = 1,2,\ldots,n$, where the parameters of  \eqref{eq:kc}, \eqref{eq:kcbc} 
are assumed to be continuous in $y$ and $\eta$ so that the pointwise evaluations are 
well-defined. Under \eqref{eq:afn1}, we have shown in \cite{HumBek24arxiv, 
HumBekCDC24} that the solution to 
\eqref{eq:kn}, \eqref{eq:knbc} converges to the solution of \eqref{eq:kc}, \eqref{eq:kcbc} (in the 
$L^2$ sense in $y$; see \cite[Lem. 4.3]{HumBek24arxiv})\footnote{This relies on 
interpreting the continuum parameters of 
\eqref{eq:kc}, \eqref{eq:kcbc} as the limits of sequences of functions defined as piecewise 
constant in $y$, on intervals of the form $((i-1)/n,i/n]$ (see \cite[Lem. 
4.2]{HumBek24arxiv}).}. Consequently, the solution of 
\eqref{eq:kc}, \eqref{eq:kcbc} can be used to 
approximate the solution to \eqref{eq:kn}, \eqref{eq:knbc} and to construct an exponentially 
stabilizing feedback gain for a large-scale $n+1$ hyperbolic system of PDEs, which is our main 
motivation of studying solving \eqref{eq:kc}, \eqref{eq:kcbc}.

\section{Power Series Approach to Solving Continuum Kernel Equations} \label{sec:pow}

In this section, we present a power series-based approach to solving the continuum kernel 
equations \eqref{eq:kc}, \eqref{eq:kcbc} for analytic parameters $\mu, \lambda, \sigma, \theta, 
W$, and $q$. We show that, in this case, the solution $(k, \bar{k})$ to \eqref{eq:kc}, 
\eqref{eq:kcbc} is analytic, and propose a numerical procedure to approximate the solution to 
arbitrary accuracy. Corresponding approaches have been adopted to solve various kernel PDEs 
\cite{VazCheCDC23, VazZha23, LinVaz24arxiv}, but our developments are new for continua of 
kernel PDEs \eqref{eq:kc}, \eqref{eq:kcbc}. Moreover, we provide a computational metric to 
assess the accuracy of the approximate solution, and consider computing reduced-order 
approximations to reduce the computational complexity of the proposed approach.

\subsection{Convergence of Power Series Representation}

As per \cite{VazCheCDC23}, the idea is to find the solution to
\eqref{eq:kc}, \eqref{eq:kcbc} as a power series, which in
case of $k(x,\xi,y)$ is a triple power series 
\begin{equation}
  \label{eq:kps}
  k(x,\xi,y) = \sum_{\ell=0}^{\infty}\sum_{i=0}^{\infty}\sum_{j=0}^iK_{ij\ell}x^{i-j}\xi^jy^{\ell},
\end{equation}
whereas for $\bar{k}(x,\xi)$ the power series representation is
\begin{equation}
  \label{eq:kbps}
  \bar{k}(x,\xi) = \sum_{i=0}^{\infty}\sum_{j=0}^i\bar{K}_{ij}x^{i-j}\xi^j.
\end{equation}
 Similarly, the parameters of \eqref{eq:kc}, 
\eqref{eq:kcbc} are represented by the series
\begin{subequations}
\label{eq:ppc}%
\begin{align}
  \lambda(x,y) & = \sum_{i=0}^\infty\sum_{j=0}^i
	\lambda_{ij}x^{i-j}y^j, \\
	\mu(x) & = \sum_{i=0}^\infty \mu_ix^i, \\
\theta(\xi,y) & = \sum_{i=0}^{\infty}\sum_{j=0}^i
  \theta_{ij}\xi^{i-j}y^j, \\
W(\xi,y) & = \sum_{i=0}^{\infty}\sum_{j=0}^i
           W_{ij}\xi^{i-j}y^j, \\
\sigma(x,\eta,y) & =
                   \sum_{\ell=0}^{\infty}\sum_{i=0}^{\infty}\sum_{j=0}^i\sigma_{ij\ell}x^{i-j}\eta^jy^{\ell},
  \\
  q(y) & = \sum_{i=1}^{\infty}q_iy^i,
\end{align}
\end{subequations}
where the coefficients are obtained from the Taylor series of the
respective parameters. Similarly to \cite{VazCheCDC23}, we consider the parameters and the 
kernels appearing in \eqref{eq:kps}--\eqref{eq:ppc} complex-valued, and by a polydisk we refer to 
$\mathcal{D}_L \times \mathcal{D}_L \times \mathcal{D}_L$ (or $\mathcal{D}_L \times 
\mathcal{D}_L$ if only two spatial variables are involved), where $\mathcal{D}_L$ is a 
complex-valued open disk centered at the origin of radius $L$, i.e.,
\begin{equation}
\label{eq:DL} 
\mathcal{D}_L = \left\{z \in \mathbb{C}: |z| < L\right\}.
\end{equation}
For the power series in \eqref{eq:kps}--\eqref{eq:ppc} to converge in the domain of definition 
of the kernel equations \eqref{eq:kc}, 
the involved functions have to be analytic on polydisks with radius larger than 
one.\footnote{Alternatively, the power series can be developed with 
	respect to some other point than the origin, in which case polydisks with smaller radii are 
	sufficient, provided that the prism $\mathcal{P}$ (or the triangle $\mathcal{T}$ in 
	case of $\bar{k}$) is covered.}

By \cite[Thm 3]{AllKrs23arxiv}, under sufficient continuity assumptions on the parameters of 
\eqref{eq:kc}, \eqref{eq:kcbc}, there exists a unique solution $(k,\bar{k})$ to 
\eqref{eq:kc}, \eqref{eq:kcbc}. However, representing the solution as a power series requires the 
solution to be analytic, so that stronger assumptions have to be imposed on the parameters of 
\eqref{eq:kc}, \eqref{eq:kcbc}. On the other hand, if an analytic solution exist, it is uniquely given 
by the power series \eqref{eq:kps}, \eqref{eq:kbps}, because of the uniqueness of the solution 
and the Taylor series representation. Similar to the results of \cite{VazCheCDC23}, 
we utilize the well-posedness result \cite[Thm 3]{AllKrs23arxiv} of \eqref{eq:kc}, \eqref{eq:kcbc} 
to state the following.

\begin{theorem}
	\label{thm:pswp}
	If the parameters of \eqref{eq:kc}, \eqref{eq:kcbc} are analytic on polydisks with radii larger 
	than one, so that they can be represented as the series in \eqref{eq:ppc}, and $|\lambda(x,y)| > 
	0$ for all $x,y \in [0,1]$, $|\mu(x)| > 0$ for all $x \in [0,1]$, then the series 
	defined in \eqref{eq:kps}, \eqref{eq:kbps} converge. That is, they define analytic functions on 
	polydisks with radii larger than one, which are the unique solution to the kernel equations 
	\eqref{eq:kc}, \eqref{eq:kcbc}.
	\begin{proof}
		The proof follows the same steps as the argument for the analogous result for the $1+ 
		1$ system in \cite[Thm 3]{VazCheCDC23}, that is, by complexifying the kernel 
		well-posedness proof of \cite[Thm 3]{AllKrs23arxiv}. In more detail, extend the successive 
		approximation series given in \cite[(109), (110)]{AllKrs23arxiv} to polydisks, which 
		requires 
		considering the integrals along the characteristic curves as line integrals in complex 
		spaces.\footnote{The complexified equations are identical to \cite[(109), 
		(110)]{AllKrs23arxiv}; only the interpretation (real vs. complex) changes.} 
		As all the parameters are 
		assumed to be analytic, the products and inner products appearing inside the
		integrals are analytic as well, and such integrals of analytic functions are independent of the 
		integration path. Consequently, it follows by recursion that each term in the successive 
		approximation series is analytic, being composed of integrals and (inner) products of analytic 
		functions. The uniform convergence of the complexified series of successive approximations 
		follows by the Weierstrass M-test (see, e.g., \cite[Thm 7.10]{RudBook76}) after derivation of 
		similar estimates as in 
		\cite[(111)--(116), (121)--(128)]{AllKrs23arxiv} for the complexified parameters.
	\end{proof}
\end{theorem}

\subsection{Computation of Approximate Solution with Truncated Power Series} 
\label{pow:alg}

Based on Theorem~\ref{thm:pswp}, we can construct arbitrarily accurate approximations of the 
solution to \eqref{eq:kc}, \eqref{eq:kcbc}, assuming that the involved parameters are analytic, by 
truncating the infinite series  \eqref{eq:kps}--\eqref{eq:ppc} up to some sufficiently high order 
$N$. This results in finding an approximate solution $(k^\mathrm{a}, 
\bar{k}^\mathrm{a})$ to 
\eqref{eq:kc}, \eqref{eq:kcbc} as
\begin{subequations}
	\label{eq:solpsa}
  \begin{align}
  k^\mathrm{a}(x,\xi,y) & = \sum_{\ell=0}^N\sum_{i=0}^{N-\ell}\sum_{j=0}^i 
  K^\mathrm{a}_{ij\ell}x^{i-j}\xi^jy^{\ell},   
  \label{eq:kpsa} \\
	\bar{k}^\mathrm{a}(x,\xi) &  = \sum_{i=0}^N\sum_{j=0}^i \bar{K}^\mathrm{a}_{ij}x^{i-j}\xi^j, 	
	\label{eq:kbpsa}
	\end{align}
\end{subequations}
for approximate parameters
\begin{subequations}
	\label{eq:ppca}%
	\begin{align}
		\lambda^\mathrm{a}(x,y) & = \sum_{i=0}^N\sum_{j=0}^i
		\lambda_{ij}x^{i-j}y^j, \\
		\mu^\mathrm{a}(x) & = \sum_{i=0}^N \mu_ix^i, \\
		\theta^\mathrm{a}(\xi,y) & = \sum_{i=0}^N\sum_{j=0}^i 
		\theta_{ij}\xi^{i-j}y^j, \\
		W^\mathrm{a}(\xi,y) & = \sum_{i=0}^N\sum_{j=0}^i
		W_{ij}\xi^{i-j}y^j, \\
		\sigma^\mathrm{a}(x,\eta,y) & =
		\sum_{\ell=0}^N\sum_{i=0}^{N-\ell}\sum_{j=0}^i\sigma_{ij\ell}x^{i-j}\eta^jy^{\ell},
		\\
		q^\mathrm{a}(y) & =  \sum_{i=1}^Nq_iy^i.
	\end{align}
\end{subequations}
Note that while the parameters are approximated by truncating the series \eqref{eq:ppc}, the 
coefficients of \eqref{eq:solpsa} may not (for any finite $N$) match the corresponding 
coefficients of \eqref{eq:kps}, \eqref{eq:kbps} (which are unknown). A numerical procedure to 
compute the coefficients of the 
truncated power series approximation \eqref{eq:solpsa} by least-squares is presented in 
Algorithm~\ref{alg:ps}, and the asymptotic convergence of the proposed approximation to the 
exact solution \eqref{eq:kps}, \eqref{eq:kbps} is established in 
Proposition~\ref{thm:alg}.\footnote{As presented here, the algorithm can be 
implemented, e.g., in MATLAB by utilizing symbolic computations, and such an 
implementation is available at 
\url{https://github.com/jphumaloja/Continuum-Kernels-Power-Series/}. Another approach would 
be to 
adapt the double series vector-matrix framework from \cite[Sect. II.B]{LinVaz24arxiv}, which is 
computationally more efficient, but potentially also more laborious to implement for the triple 
power series appearing in our computations.}
\begin{algorithm}
	\DontPrintSemicolon
	\KwData{Parameters $\lambda,\mu,\sigma,\theta,W,q$ of \eqref{eq:kc}, \eqref{eq:kcbc}.}
	\KwResult{Truncated power series approximation \eqref{eq:solpsa} for $(k,\bar{k})$.}
	\Init{}{Choose approximation order $N$.\;
	Construct the truncated series \eqref{eq:solpsa}, \eqref{eq:ppca} (with unknowns the 
	coefficients of \eqref{eq:solpsa}).\;
	Insert the truncated series into \eqref{eq:kc}, \eqref{eq:kcbc}.\;
	Store the coefficients of each term $x^i\xi^{j-i}y^\ell$ from \eqref{eq:kc}, \eqref{eq:kcbc} into a 
	list.\;
	Initialize matrix $\mathbf{A}$ and vector $\mathbf{b}$ based on the number of unknown 
	coefficients $K^\mathrm{a}_{ij\ell},\bar{K}^\mathrm{a}_{ij}$ and the length of the list.
	}
	\While{not at the end of the list}{
		Check next element in the list of coefficients.\;
		\eIf{coefficient has $K^\mathrm{a}_{ij\ell}$ or $\bar{K}^\mathrm{a}_{ij}$}{
			Insert the numerical value to the corresponding position in $\mathbf{A}$.\;
		}{
			Insert the numerical value to the corresponding position in $\mathbf{b}$.\;
		}
	}
	Solve for the unknown coefficients from $\mathbf{Ax}^\mathrm{a}=\mathbf{b}$ by least 
	squares.\; 
	Insert the obtained values $\mathbf{x}^\mathrm{a}$ to \eqref{eq:solpsa}.\;
	\label{alg:ps}
	\caption{Computation of the power series approximation \eqref{eq:solpsa}.}
\end{algorithm}
\begin{proposition}
	\label{thm:alg}
	Under the conditions of Theorem \ref{thm:pswp}, the approximate solution \eqref{eq:solpsa} 
	produced by Algorithm~\ref{alg:ps} converges asymptotically to the exact solution 
	\eqref{eq:kps}, \eqref{eq:kbps} as $N\to\infty$.
	\begin{proof}
	For any $N$, refer to the linear system of equations of Algorithm~\ref{alg:ps} as 
	$\mathbf{A}_N\mathbf{x}_N = \mathbf{b}_N$, and denote the coefficients of the (exact) 
	truncated power series \eqref{eq:kps}, \eqref{eq:kbps} of order $N$ by 
	$\mathbf{x}^\mathrm{e}_N$. 
	By 	definition, a least-squares  solution to $\mathbf{A}_N\mathbf{x}_N = \mathbf{b}_N$, 
	denoted by $\mathbf{x}_N^\mathrm{a}$, minimizes $\|\mathbf{A}_N\mathbf{x}_N - 
	\mathbf{b}_N\|_2$. Thus, for all $N$, we have 
	\begin{equation}
		\label{eq:algthm1}
		\| \mathbf{A}_N\mathbf{x}_N^\mathrm{a} - \mathbf{b}_N\|_2 \leq \| 
		\mathbf{A}_N\mathbf{x}_N^\mathrm{e} - 
		\mathbf{b}_N\|_2.
	\end{equation}
	In the limit case $N\to\infty$, by the uniqueness of the solution to \eqref{eq:kc}, 
	\eqref{eq:kcbc} 
	and the uniqueness of the power series representation, the coefficients 
	$\mathbf{x}_\infty^\mathrm{e}$ of 
	\eqref{eq:kps}, \eqref{eq:kbps} uniquely solve the (infinite) set of linear equations 
	$\mathbf{A}_\infty \mathbf{x}_\infty = \mathbf{b}_\infty$, i.e., $\|\mathbf{A}_\infty 
	\mathbf{x}_\infty^\mathrm{e} - \mathbf{b}_\infty\|_{2} = 0$. Combining this with 
	\eqref{eq:algthm1} gives 
	$\mathbf{x}_\infty^\mathrm{a} = \mathbf{x}_\infty^\mathrm{e}$. Thus, the coefficients of 
	\eqref{eq:solpsa} tend 
	to the coefficients of \eqref{eq:kps}, \eqref{eq:kbps} as $N\to\infty$, and the claim follows.
\end{proof}
\end{proposition}

As stated in Algorithm~\ref{alg:ps}, the solution to the kernel equations is obtained by 
replacing the parameters and the solution in the kernel equations \eqref{eq:kc}, \eqref{eq:kcbc} 
by the power series approximations \eqref{eq:solpsa}, \eqref{eq:ppca}. In turn, solving for the 
unknown coefficients $K^\mathrm{a}_{ij\ell}$ and $\bar{K}^\mathrm{a}_{ij}$ such that 
the equations are 
satisfied, matching the coefficients of powers of $x$, $\xi$, and $y$ appearing in the
series. In particular, as also noted in \cite[Sect. III.D]{VazCheCDC23}, it is reasonable to express 
the boundary condition \eqref{eq:kcbca} as
\begin{equation}
	\label{eq:psbc1}
	(\lambda(x,y) + \mu(x))k(x,x,y) = -\theta(x,y),
\end{equation}
so that both sides are polynomials in $x$ and $y$, because the original form \eqref{eq:kcbca} is 
not directly compatible with the power series approach.\footnote{Alternatively, one could 
compute a power series approximation separately for the right-hand side of \eqref{eq:kcbca}. 
However, one should note that the analyticity of $\lambda, \mu$, and $\theta$ does not generally 
imply that the fraction on the right-hand side of \eqref{eq:kcbca} would be analytic. Hence, the 
form \eqref{eq:psbc1} is preferable, where the functions are replaced by their power series 
approximation when solving for the unknown coefficients.} Finally, we note that the set of 
linear equations considered in Algorithm~\ref{alg:ps} is usually over-determined. This is due to the 
parameters and solution being approximated by power series of order~$N$, so that any 
multiplication  in \eqref{eq:kc}, \eqref{eq:kcbc} results in the appearance of terms up to order 
$2N$, which contribute to the equations that have to be satisfied by the series coefficients.

By virtue of Proposition~\ref{thm:alg}, the accuracy of the approximate solution obtained with 
Algorithm~\ref{alg:ps} can be assessed based on the residual $\|\mathbf{Ax}-\mathbf{b}\|_2$ of 
the least-squares fit. If the residual is large, then one should increase the approximation order 
$N$. Finally, as the convergence of power series is uniform (in space), the obtained approximate 
solution is arbitrarily close to the exact solution pointwise, meaning that we also get arbitrarily 
accurate approximations for the state feedback gains $k(1,\xi,y)$ and $\bar{k}(1,\xi)$ employed in 
the backstepping control law (see \eqref{eq:Uc}). The spatial error of the approximation can 
in fact be estimated, 
similarly to the remainder of the truncated series \eqref{eq:kps}, \eqref{eq:kbps}, by expressions 
of the form (see \cite[Prop. 1.2.2]{Leb23Book})
\begin{equation}
	|R_N(x,\xi,y)| \leq M \frac{L^{D(N+1)}}{(N+1)!},
\end{equation} 
where $D = 2,3$ for \eqref{eq:kbps}, \eqref{eq:kps}, respectively, and $M$ is a bound depending 
on the derivatives of order $N+1$ of the approximated function. 

The bound $M$ can 
(technically) be estimated by deriving the kernel equations for the derivatives of order $N+1$ of 
$(k,\bar{k})$ by differentiating \eqref{eq:kc}, \eqref{eq:kcbc} with respect to $(x,\xi,y)$ (see, e.g., 
\cite[App. A.4]{CorVaz13} for the case of $2\times  2$ kernels). The resulting kernel equations are 
essentially of the same form as \eqref{eq:kc}, \eqref{eq:kcbc} (with the lower-order derivatives of 
$(k,\bar{k})$ acting as additional parameters), and thus, their solution can be upper-bounded 
using the method of successive approximations as in \cite[Sect. VI]{AllKrs23arxiv}.  However, in 
order to compute such a bound, all derivatives of $(k,\bar{k})$ up to order $N+1$ would need to 
be estimated, implying that $\frac{(N+2)(N+3)}{2}$ estimates for the derivatives of $\bar{k}$ and 
$\frac{(N+1)(N+2)(2N+12)}{12}+N+2$ estimates for the 
derivatives of $k$ would have to be computed, which has the same computational complexity 
$\mathcal{O}(N^3)$ as the 
proposed power series procedure for actually computing the kernels themselves. (see 
Section~\ref{pow:cc}). In other words, it would be
more practically useful to 
evaluate the truncation error from the least-squares error of Algorithm~\ref{alg:ps}, rather than 
trying to quantify it through the power series residual.

\subsection{Computational Complexity of the Power Series-Based Kernels Approximation} 
\label{pow:cc}

The number of unknown coefficients in the truncated power series \eqref{eq:solpsa} is
$\displaystyle \frac{(N+1)(N+2)}{2}$ for $\bar{K}^\mathrm{a}_{ij}$ and 
\begin{equation} 
  \sum_{i=0}^N \frac{(i+1)(i+2)}{2}
    = \frac{N(N+1)(2N+10)}{12} + N+1,
\end{equation}
for $K^\mathrm{a}_{ij\ell}$, meaning that we have that many unknowns to solve for, whose 
number 
increases proportionally to $N^3$. Thus, the computational complexity of the triple power series 
approach grows cubically with the approximation order $N$. This is consistent with quadratic 
growth reported in \cite{VazCheCDC23} for the double power series.

If we view \eqref{eq:kc}, \eqref{eq:kcbc} as an approximation of \eqref{eq:kn}, \eqref{eq:knbc} and 
think of solving \eqref{eq:kn}, \eqref{eq:knbc} using a power series approach, this would involve 
$n+1$ double power series, one for each of the components $k^i$ (similar to $\bar{k}$ in the 
continuum case), meaning that the power series approach to solve \eqref{eq:kn}, \eqref{eq:knbc} 
would involve solving for $\displaystyle (n+1)\frac{(N+1)(N+2)}{2}$ unknown coefficients. Thus, 
the growth in computational complexity to solve \eqref{eq:kn}, \eqref{eq:knbc} based on a power 
series approach is  $\mathcal{O}(nN^2)$, whereas for solving the continuum kernel equations 
\eqref{eq:kc}, \eqref{eq:kcbc} it is $\mathcal{O}(N^3)$. Hence, whenever it is possible to choose 
(roughly) $N < n$ without compromising the approximation accuracy, the continuum 
approximation approach of solving \eqref{eq:kc}, \eqref{eq:kcbc} reduces the computational 
complexity of the power series approach compared to solving \eqref{eq:kn}, \eqref{eq:knbc}. 
More importantly, in the case of solving \eqref{eq:kc}, \eqref{eq:kcbc}, the computational 
complexity does not scale with the number of components of the $n+1$ system.

In the next subsection, we discuss a potential order 
reduction method with respect to the ensemble variable $y$, which may reduce the complexity of 
the triple power series from $\mathcal{O}(N^3)$ to $\mathcal{O}(N^2)$. That is, instead of a 
generic order for the truncated power series \eqref{eq:solpsa}, we employ a separate, lower order 
$N_y<N$ for the powers of $y$, which leads to a smaller number of unknown coefficients, and 
hence, to
computational complexity of order $\mathcal{O}(N_yN^2)$. Thus, assuming that the order $N_y$ 
can be kept constant, the computational complexity only grows quadratically with respect to $N$ 
as opposed to cubically in \eqref{eq:solpsa}. Moreover, compared to the complexity 
$\mathcal{O}(nN^2)$ of solving \eqref{eq:kn}, \eqref{eq:knbc} with a power series approach, the 
reduced-order approach with respect to $y$ further reduces computational complexity whenever 
(roughly) $N_y < n$. In both cases of computing the power series corresponding to \eqref{eq:kc}, 
\eqref{eq:kcbc}, computational complexity does not grow with $n$, in contrast to the case of 
computing the power series corresponding to \eqref{eq:kn}, \eqref{eq:knbc}. The computational 
complexities of the three considered approaches are 
summarized in Table~\ref{tab:ccc} below.
\begin{table}[!htb]
	\begin{center}
		\begin{tabular}{c | c c c} 
			Approach & Full \eqref{eq:kc}, \eqref{eq:kcbc} & \eqref{eq:kn}, \eqref{eq:knbc} & R-O 
			$y$ \eqref{eq:kc}, \eqref{eq:kcbc} 	\\ 
			\hline \\ [-11pt]
			Complexity & $\mathcal{O}(N^3)$ & $\mathcal{O}(nN^2)$ & $\mathcal{O}(N_yN^2)$ 
		\end{tabular}
	\end{center}
	\caption{Computational complexities of the power series approaches to solving \eqref{eq:kn}, 
	\eqref{eq:knbc} and \eqref{eq:kc}, \eqref{eq:kcbc} with a potentially reduced-order (R-O) in 
	$y$ series.}
	\label{tab:ccc}
\end{table}

As our continuum approximation approach asks 
	for sufficiently large $n$, we, in general, expect $N$ or $N_y$ to be much smaller than $n$, or, 
	more accurately, our approach, theoretically, relies on the assumption that $n$ is (sufficiently) 
	large. Thus, a more fair comparison would be to compare $N$ with $N_y$ so that, in addition to 
	having computational complexity not growing with $n$, we further reduce computational 
	complexity from $N^3$ to $N_y N^2$, despite that, from a purely algebraic viewpoint, we have 
	to also compare $N_y$ with $n$ so that to illustrate (algebraically) that complexity $O(n N^2)$ 
	is worse also as compared with $O(N_y N^2)$.

\subsection{On Order Reduction with Respect to $y$} \label{pow:red}

Instead of having a generic degree 
$N$ for the power series approximations \eqref{eq:solpsa}, we can choose a lower degree for 
some of the spatial 
variables. Considering that $y$ is the ensemble variable (when interpreting \eqref{eq:kc}, 
\eqref{eq:kcbc} as the continuum of \eqref{eq:kn}, \eqref{eq:knbc}), and there 
is some freedom involved in constructing the continuum parameters satisfying \eqref{eq:afn1} 
(see \cite{HumBek24arxiv}), it 
may be possible to obtain sufficiently accurate power series approximations with a much lower 
order in $y$. One such particular case is when all the parameters are, or can be approximated 
by, low-order polynomials in $y$.
Based on the power series approximation in \eqref{eq:kpsa}, we can construct a 
reduced-order approximation in $y$ as 
\begin{equation}
	\label{eq:kpsay}
	k^{\mathrm{a}_y}(x,\xi,y) = \sum_{\ell=0}^{N_y}\sum_{i=0}^{N-\ell}\sum_{j=0}^i 
	K^{\mathrm{a}_y}_{ij\ell}x^{i-j}\xi^jy^{\ell},
\end{equation}
for some $N_y < N$. 
As reported in Table~\ref{tab:ccc}, the respective computational complexity to solve for the 
unknown coefficients is $\mathcal{O}(N_yN^2)$.

Assuming that the parameters $\lambda, \theta$, and $\sigma$ are, or can be approximated by,  
low-order polynomials in $y$, one can try to approximate the continuum kernels with 
reduced-order $N_y$ 
in $y$. The reason as to why a lower order $N_y$ in $y$, than in the other 
spatial variables, would still provide an accurate approximation, is that \eqref{eq:kc}, 
\eqref{eq:kcbc} is not a differential equation in $y$. Hence, an 
order in $y$ similar to the respective order of the parameters, determining the spatial 
dependence of \eqref{eq:kca} in $y$, may be enough to approximate the solution to  
\eqref{eq:kc}, \eqref{eq:kcbc} accurately. On the other hand, considering that particularly
\eqref{eq:psbc1} has to be satisfied, the order of $\lambda(x,y)k(x,x,y)$ has to be at least the 
same as the order of $\theta(x,y)$ in $y$. That is, if the orders of $\lambda(x,y)$ and 
$\theta(x,y)$ in $y$ are $N_\lambda$ and $N_\theta$, respectively, then necessarily $N_y \geq 
N_\theta - N_\lambda$, which provides a lower bound for $N_y$.\footnote{We have observed in 
numerical experiments that, in general, this 
lower bound may be sufficient particularly when $\lambda$ 
is constant in $y$, i.e., $N_\lambda = 0$. However, for $y$-varying $\lambda$, it may be more 
difficult to determine, a priori, a sufficient approximation order $N_y$ in $y$ (without 
compromising the approximation accuracy).}

If the continuum kernel equations 
\eqref{eq:kc}, \eqref{eq:kcbc} are treated as an approximation of the large-scale kernel equations 
\eqref{eq:kn}, \eqref{eq:knbc}, a potential additional benefit of the continuum approximation can 
be gained 
particularly  when there exists continuum parameters that are low-order 
polynomials in $y$, such that the relations \eqref{eq:afn1} are satisfied for the corresponding 
large-scale parameters, which potentially results in $N_y <n$.
For example, this is the case if the large-scale parameters are parametrized by low-order 
polynomials of $i/n$, i.e., they form sequences of polynomials in $i$, in which case the continuum 
parameters satisfying \eqref{eq:afn1} can be constructed by replacing $i/n$ with $y$ in the 
expressions of the large-scale parameters.\footnote{Note that there is quite a bit of 
freedom 
in constructing the continuum approximation, stemming from the fact that the only restriction in 
obtaining the continuum parameters from the parameters of the $n+1$ system is \eqref{eq:afn1} 
(see also 
\cite{HumBek24arxiv,HumBekCDC24} for details).} We demonstrate this later on in 
Section~\ref{sec:ex2}.

In general, we note that  
\eqref{eq:afn1} is not strictly necessary for the continuum approximation to be sufficiently 
accurate (cf. \cite[Rem. 4.4]{HumBek24arxiv}). Essentially, any polynomial that approximates 
the distribution of the large-scale parameters as an ensemble (in the $L^2$ sense in $y$) is 
sufficient.  For finding such polynomial approximations in practice, many computational software 
have routines for such purposes, such as, e.g., \verb+polyfit+ in MATLAB, which can be 
potentially utilized in constructing such polynomial continuum approximations. We demonstrate 
this approach as well in Section~\ref{sec:ex2}.

\section{On Explicit Solution to Continuum Kernel Equations} \label{sec:exp}

In this section, we identify sufficient conditions on the parameters of the kernel equations 
\eqref{eq:kn}, 
\eqref{eq:knbc}, under which the continuum kernel equations \eqref{eq:kc}, 
\eqref{eq:kcbc} can be solved explicitly. The main idea towards this is to look for separable 
solutions to the continuum kernel equations, the existence of which requires certain assumptions 
on the continuum parameters, and respectively, on the parameters of the original large-scale 
kernel equations.\footnote{We note that these are merely sufficient conditions, so that it may 
	be possible to derive closed-form solutions under different assumptions as well. In particular, 
	as opposed to Assumption~\ref{ass:n+1}, one may consider spatially-varying $\lambda_i, 
	i=1,\ldots,n$ and 
	$\mu$. However, deriving a closed-form solution under such conditions likely requires imposing 
	additional assumptions on the parameters, as opposed to the three conditions 
	\eqref{eq:cx}--\eqref{eq:fd} we impose in Proposition~\ref{prop:esol}. Regardless, the 
	case of constant 
	transport speeds $\lambda_i, \mu$ appears at least in the examples considered in, e.g.,  
	\cite[Sect. 
	VII]{AllKrs23arxiv} and \cite[Sect. VI]{DiMVaz13}, while such an assumption holds for linearized 
	hyperbolic systems, describing, e.g., traffic or shallow water flows, around a constant 
	equilibrium point; see, e.g., \cite{BasCorBook, YuHKrs21}.} In Section~\ref{sec:ex1}, we 
	utilize the explicit solutions to benchmark the approximate power series approach presented in 
	Section~\ref{sec:pow}.

\begin{assumption}\label{ass:n+1}
	For the parameters of the kernel equations \eqref{eq:kn}, \eqref{eq:knbc}, we assume that 
	$\mu(x)  = \mu$ and $\lambda_i(x) = \lambda_i$, where $\mu$ and $\lambda_i$, for all $i 
	= 1,\ldots,n$, are positive constants. Moreover, the parameters 
	$W_i(x), \sigma_{i,j}(x), \theta_i(x)$ are of the form 
	\begin{subequations}
		\label{eq:wstncond}%
		\begin{align}
			W_i(x) & =  w_iW_x(x), \\
			\sigma_{i,j}(x) & = s_{1,i}s_{2,j}\sigma_{x}(x), \\
			\theta_i(x) & = \vartheta_i\theta_x(x),
		\end{align}
	\end{subequations}
	for some constants $w_i, s_{1,i}, s_{2,j}, \vartheta_i$ for all $i,j = 1,\ldots,n$.
\end{assumption}
Under Assumption~\ref{ass:n+1}, the parameters of the continuum kernel equations \eqref{eq:kc}, 
\eqref{eq:kcbc} can be constructed to be separable (with $\lambda$ constant in $x$), i.e., 
there exist functions $\lambda_y, 
W_y, \theta_y,\sigma_y,\sigma_{\eta}$ such that
\begin{subequations}
	\label{eq:csep}%
	\begin{align}
		\lambda(x,y) & = \lambda_y(y), \\
		W(x,y) & =  W_x(x)W_y(y), \\
		\sigma(x,y,\eta) & = \sigma_{x}(x)\sigma_y(y)\sigma_{\eta}(\eta), \\
		\theta(x,y) & = \theta_x(x)\theta_y(y),
	\end{align}
\end{subequations}  
where
\begin{subequations}
	\label{eq:afn}%
	\begin{align}
		\lambda_y(i/n) & = \lambda_i, \\
		W_y(i/n) & = w_i, \\
		\theta_y(i/n) & = \vartheta_i, \\
		\sigma_y(i/n) & = s_{1,i}, \\	
		\sigma_{\eta}(i/n) & = s_{2,i},
	\end{align}%
\end{subequations}
for all $i=1,\ldots,n$. Now, we can identify additional conditions that have to be satisfied by 
the parameters appearing in \eqref{eq:csep}, in order to guarantee the existence of a separable, 
closed-form solution to the continuum kernel 
equations \eqref{eq:kc}, \eqref{eq:kcbc}. This is  formulated in the following proposition. The 
motivation for the proposition is that we can impose corresponding conditions on the parameters 
of the large-scale kernel equations \eqref{eq:kn}, \eqref{eq:knbc}, so that the closed-form 
continuum solution to \eqref{eq:kc}, \eqref{eq:kcbc} provides an approximation for the solution to  
\eqref{eq:kn}, \eqref{eq:knbc}. This is discussed in detail after Proposition~\ref{prop:esol} in 
Remark~\ref{rem:n1}.

\begin{proposition}
	\label{prop:esol}
	Let Assumption~\ref{ass:n+1} hold  so that $\mu$ is constant and the 
	continuum parameters are separable as in \eqref{eq:csep}, satisfying \eqref{eq:afn}.
	Additionally assume that there exists a constant $c_x$ satisfying
\begin{align}
\label{eq:cx}
c_x & = 
\mu\sigma_x(0)\frac{\sigma_\eta(y)}{\theta_y(y)}\int\limits_0^1 
\frac{\sigma_y(\eta)\theta_y(\eta)}{\lambda(\eta) + \mu}d\eta + 
\frac{\mu\lambda(y)}{\lambda(y)+\mu}
\frac{\theta'_x(0)}{\theta_x(0)} \nonumber \\ & \qquad +
\theta_{x}(0)\int\limits_0^1\lambda(\eta)q(\eta)\frac{\theta_y(\eta)}{\lambda(\eta) + 
	\mu}d\eta,
\end{align}
independently of $y$, and a related function $f$ satisfying 
\begin{align}
	\label{eq:f}
	f(\xi) & = \sigma_x(\xi)\frac{\sigma_\eta(y)}{\theta_y(y)}\int\limits_0^1 
	\frac{\sigma_y(\eta)\theta_y(\eta)}{\lambda(\eta) + \mu}d\eta \nonumber\\ & 
	\qquad  -\frac{c_x}{\mu} + \frac{\lambda(y)}{\lambda(y)+\mu}\frac{\theta_x'(\xi)}{\theta_x(\xi)},
\end{align}
for all $\xi \in [0,1]$, independently of $y$, and 
\begin{equation}
	\label{eq:fd}
	f'(\xi) = W_x(\xi)\theta_x(\xi)\int\limits_0^1\frac{W_y(y)\theta_y(y)}{\lambda(y) + \mu}dy,
\end{equation}
for all $\xi \in [0,1]$. Then, the solution to the continuum kernel equations \eqref{eq:kc}, 
\eqref{eq:kcbc} is given by
\begin{subequations}
	\label{eq:ksep}
	\begin{align}
		k(x,\xi,y) & = -\exp \left(\frac{c_x}{\mu}x \right)\exp
		\left(-\frac{c_x}{\mu}\xi \right)\theta_x(\xi)\frac{\theta_y( y)}{\lambda(y)+\mu}, 
		\label{eq:ksepa} \\
		\bar{k}(x,\xi) & = \exp \left(\frac{c_x}{\mu}x \right)f(\xi)\exp\left(-\frac{c_x}{\mu}\xi \right). 
		\label{eq:ksepb}
	\end{align}
\end{subequations}
	\begin{proof}
		The proof can be found in \ref{app:esol}.		
	\end{proof}
\end{proposition}

\begin{remark}
\label{rem:prop}
If $\lambda$ is assumed to be constant, it is possible to make the somewhat implicit conditions of 
Proposition~\ref{prop:esol} explicit. In this case, the existence of the constant $c_x$ and the 
function $f$ reduces to the existence of a constant $c_y$ satisfying
\begin{equation}
	\label{eq:cy}
	c_y = \frac{\sigma_\eta(y)}{\theta_y(y)}\int\limits_0^1\sigma_{y}(\eta)\theta_y(\eta)d\eta,
\end{equation}
independently of $y$, i.e., either there exists a constant $c$ such that $\sigma_\eta = c\theta_y$ 
or the integral appearing in \eqref{eq:cy} is zero. If such $c_y$ exists, the constant $c_x$ is then 
given by
\begin{equation}
	\resizebox{.99\columnwidth}{!}{$
		\label{eq:cx2}
		\displaystyle c_x = \frac{\mu}{\lambda+\mu}\left( c_y\sigma_x(0) + \lambda
		\frac{\theta'_x(0)}{\theta_x(0)} +
		\frac{\lambda}{\mu}\theta_{x}(0)\int\limits_0^1q(y)\theta_y(y)dy
		\right),
		$}
\end{equation}
and the function $f$ is given by
\begin{equation}
\label{eq:f2}
f(\xi) = \frac{c_y}{\lambda+\mu}\sigma_x(\xi) - \frac{c_x}{\mu} +
\frac{\lambda}{\lambda+\mu}\frac{\theta_x'(\xi)}{\theta_x(\xi)}.
\end{equation}
Thus, the condition \eqref{eq:fd} reduces to
\begin{align}
	\label{eq:sepcond}
	c_y\sigma_x'(\xi) +
	\lambda\frac{\theta_x''(\xi)\theta_x(\xi)-\theta_x'(\xi)^2}
	{\theta_x(\xi)^2}  & = \nonumber \\
	W_x(\xi)\theta_x(\xi)	\int\limits_0^1W_y(y)\theta_y(y)dy, &
\end{align}
for all $\xi\in [0,1]$.
\end{remark}

\begin{remark} \label{rem:n1}
	The special case of constant $\lambda$ in the continuum translates to every $\lambda_i$ 
	being the same for all $i=1,\ldots,n$. Similarly, the other assumptions considered in 
	Remark~\ref{rem:prop} can be translated 
	to conditions on the parameters of the large-scale equations shown in \eqref{eq:wstncond}. 
	Regarding \eqref{eq:cy}, we either need that there exists a constant $c$ such that
	$\vartheta_i = cs_{2,i}$ for all $i =1,2,\ldots,n$, or 
	\begin{equation}
		\label{eq:cy0d}
		\lim_{n\to\infty} \frac{1}{n}\sum_{i=1}^n s_{1,i}\vartheta_i = 0,
	\end{equation} 
	where the limit converges to the integral\footnote{Formally, to the Riemann integral, and due to 
		continuity of the parameters, to the standard Lebesgue integral.}  in \eqref{eq:cy}, 
	by construction of the continuum parameters. Similarly, the condition 
	\eqref{eq:sepcond} can be written in terms of \eqref{eq:wstncond} by replacing the integral 
	with the corresponding infinite sum, i.e.,
	\begin{align}
		\label{eq:wstncond2}
		c_y\sigma_x'(\xi) +
		\lambda\frac{\theta_x''(\xi)\theta_x(\xi)-\theta_x'(\xi)^2}
		{\theta_x(\xi)^2}  & = \nonumber \\
		W_x(\xi)\theta_x(\xi)\lim_{n\to\infty}\frac{1}{n}\sum_{i=1}^n 
		w_i\vartheta_i, &
	\end{align}
	for all $\zeta \in [0,1]$. 
\end{remark}
\begin{remark}
	Even though it is reasonable for one to think that Assumption~\ref{ass:n+1} together with  the 
	assumptions in Remark~\ref{rem:n1}
	lead to closed-form solutions for the exact kernels PDEs \eqref{eq:kn}, \eqref{eq:knbc}, 
	corresponding to the $n+1$ system, as well, this is not the case. For example, in more technical 
	terms, even when $\displaystyle \int\limits_0^1 \sigma_y(y) \theta_y(y)dy = 0$, which
	corresponds to $c_y = 0$ in \eqref{eq:cy}, this does not imply 
	that the respective sum in \eqref{eq:cy0d} would necessarily be zero for any finite $n$ (it is 
	zero, but only in the limit $n\to\infty$). In fact, such types of continuum approximations (e.g., of 
	sums by 
	integrals) is the key to obtain a closed-form solution to the continuum kernel PDE (in contrast to 
	the exact kernels PDEs), which even though results in applying to the $n+1$ system 
	approximate control 
	kernels, it remains stabilizing in closed loop. In particular, as we have shown in 
	\cite{HumBek24arxiv,HumBekCDC24}, an approximate solution 
	$(k_\mathrm{a}^i)_{i=1}^{n+1}$ to  
	\eqref{eq:kn}, \eqref{eq:knbc} based on the 
	solution to \eqref{eq:kc}, \eqref{eq:kcbc} is given by
	\begin{subequations}
		\label{eq:kappr}
		\begin{align}
			k_\mathrm{a}^i(x,\xi) & = k(x,\xi,i/n), \qquad i = 1,2,\ldots,n, \\
			k_\mathrm{a}^{n+1}(x,\xi) & = \bar{k}(x,\xi).
		\end{align}
	\end{subequations}
\end{remark}

\begin{remark}
 The conditions in Proposition~\ref{prop:esol} or Remark~\ref{rem:prop} may be restrictive 
 (although it is 
straightforward for one to verify them, so to obtain a closed-form solution without having to 
undergo 
any power series-based or other approximation for computing the solution to \eqref{eq:kc}, 
\eqref{eq:kcbc}). However, 
whenever such an explicit solution is possible to find and, in fact, we identify a class of systems 
when this is possible via Proposition~\ref{prop:esol} and Remark~\ref{rem:prop}, it significantly 
reduces computational complexity. This is 
particularly useful in the case in which the explicit solution to the continuum kernels equations 
\eqref{eq:kc}, \eqref{eq:kcbc} is employed for computation of stabilizing kernels for the 
large-scale system \eqref{eq:n+1}, \eqref{eq:nbcuy} 
corresponding to the $n+1$ kernels equations \eqref{eq:kn}, \eqref{eq:knbc}. For example, as also 
illustrated in
\cite[Sect. V]{HumBek24arxiv} (see also \cite[Sect. V]{HumBekCDC24}) for the 
numerical example taken from \cite{AllKrs23arxiv}, it is possible to stabilize a large-scale $n+1$ 
system 
counterpart, for arbitrarily large $n$, with almost trivial computations, thanks to the availability of 
a closed-form 
solution to the continuum kernel equations (but not to the respective $n+1$ kernels equations). In 
fact, this approach (of computing stabilizing kernels for large-scale systems via continuum 
approximation) is particularly useful exactly when an explicit solution to the $n+1$ kernels 
equations \eqref{eq:kn}, \eqref{eq:knbc} is not available, but an explicit solution is available for the 
respective continuum equations \eqref{eq:kc}, \eqref{eq:kcbc} (as is the case with the example 
from \cite{AllKrs23arxiv}).

Furthermore, because the conditions in Proposition~\ref{prop:esol} concern the continuum 
parameters, there is some flexibility degree for their satisfaction in the case in which the 
continuum parameters involved (corresponding to the continuum system \eqref{eq:inf}, 
\eqref{eq:cbcuy}) are obtained as continuum 
approximation (which may not be unique) of the respective sequences of parameters of the 
large-scale system counterpart (corresponding to system \eqref{eq:n+1}, 
\eqref{eq:nbcuy}). This is discussed in detail 
in \cite[Sect. V]{HumBek24arxiv}. In particular, the parameters of \eqref{eq:kc}, 
\eqref{eq:kcbc} have to be tied to the parameters of \eqref{eq:kn}, \eqref{eq:knbc} in some way, 
which, in the present case, happens through \eqref{eq:afn1}, such that the solution to 
\eqref{eq:kc}, \eqref{eq:kcbc} can approximate the solution to \eqref{eq:kn}, \eqref{eq:knbc}. One 
instance where these conditions may be satisfied is illustrated in Example~1. Another usefulness 
of having identified a class of equations \eqref{eq:kc}, 
\eqref{eq:kcbc} that feature a closed-form solution stems from the fact that such closed-form 
solutions provide exact solutions that can be used for computation of exact errors of 
power series-based approximations, as it is done in Section~\ref{sec:ex1}.
\end{remark}

\noindent \textbf{Example 1.} The simplest case in which the conditions of 
Remark~\ref{rem:prop} are satisfied is when $c_y = 0$ and both sides of 
\eqref{eq:sepcond} are zero. Moreover, the solution 
formula becomes simpler in the particular 
case $c_x = 0$. These particular conditions (and particular value of $c_x$) are met in the 
example considered in \cite[Sect. VII]{AllKrs23arxiv}, where $\lambda = \mu = 1, W_y(y) = 
\sigma_y(y) = y-\frac{1}{2}, \theta_y(y) = y(y-1),  \theta_x(x) 
= -70\exp \left(\frac{35}{\pi^2}x\right)$, and $q(y) = \cos(2\pi y)$, in which case the solution 
\eqref{eq:ksep} becomes
\begin{subequations}
	\label{eq:exkc}%
	\begin{align}
		k(x,\xi,y) & = 35y(y-1)\exp \left( \frac{35}{\pi^2}\xi\right), \label{eq:exkca} \\
		\bar{k}(x,\xi) & = \frac{35}{2\pi^2}.
	\end{align}
\end{subequations}
The other parameter functions are $\sigma_{\eta}  = \sigma_y$, $\sigma_x(x) = x^3(x+1)$, and 
$W_x(x) = x(x+1)e^x$, although they do not appear in the expression of the solution as 
$c_x=c_y=0$ and $\displaystyle \int\limits_0^1W_y(y)\theta_y(y)dy = 0$.

\section{Numerical Experiments:  Power Series-Based Computation for Example~1} 
\label{sec:ex1}

We test the power series algorithm (full- and reduced-order in $y$) with the parameters of 
Example 1 to get an idea of the practical computational requirements and the accuracy of 
the method. The numerical experiments are performed in MATLAB 2023b on a 2023 
MacBook Pro with Apple M3 chip and 8 GB of memory. In Table~\ref{tab:numf}, numerical 
results are presented of computation of the 
full-order power 
series \eqref{eq:solpsa}, where the first column shows the order $N$ of the approximation; the 
second column shows 
the number of unknown coefficients (both $K^\mathrm{a}$ and $\bar{K}^\mathrm{a}$ 
together); the third column shows 
the number of (linear) equations (that have to be solved for determining the series coefficients); 
the fourth column shows the computational time (in seconds); the fifth column shows the residual 
of the least-squares fit when the set of linear equations is solved; and the sixth column shows the 
maximal absolute error between the computed and the exact, continuum control gains (i.e., the 
kernels \eqref{eq:exkc} evaluated at $x=1$). One can see that both the number of unknowns and 
the number of 
equations 
grow rapidly along with $N$, and hence, so does the computational time. On the flip side, the 
residual of the least-squares fit and the maximal error of the control kernel approximation become 
smaller 
as $N$ increases. However, for larger values of $N$, (notable) improvements in the 
approximation accuracy are only obtained with even $N$, e.g., the max error only reduces by 
$0.046$ when comparing $N=14$ and $N=15$, whereas for $N=15$ and $N=16$ the max error 
reduces by $0.506$. This may be related to approximating 
the even function $q(y) = \cos(2\pi y)$, as most of the other parameters (apart 
from the exponential terms in $\theta$ and $W$)  are polynomials, for which the truncated Taylor 
series are exact. 
\begin{table}[!htb]
\begin{center}
	\begin{tabular}{c | c c c c c} 
		N & \#K & \#Eq & C. time & $\|\mathbf{Ax}-\mathbf{b}\|_2$ & max error \\ 
		\hline \\ [-11pt]
		12 & $546$ & $1082$ & $11.78$ s. & $1.90$ & $16.95$ \\
		13 & $665$ & $1285$ & $13.71$ s. & $0.669$ & $8.79$ \\
		14 & $800$ & $1510$ & $18.27$ s. & $0.209$ & $0.668$ \\
		15 & $952$ & $1758$ & $25.24$ s. & $5.53\cdot10^{-2}$ & $0.622$ \\
		16 & $1122$ & $2030$ & $32.79$ s. & $1.13\cdot10^{-2}$ & $0.116$ \\
		17 & $1311$ & $2327$ & $40.90$ s. & $3.10\cdot10^{-3}$ & $9.14\cdot10^{-2}$ \\
		18 & $1520$ & $2650$ & $50.35$ s. & $6.83\cdot10^{-4}$ & $7.23\cdot10^{-3}$ \\
		19 & $1750$ & $3000$ & $62.33$ s. & $1.42\cdot10^{-4}$ & $7.61\cdot10^{-3}$ \\
		20 & $2002$ & $3378$ & $77.70$ s. & $2.82\cdot10^{-5}$ & $5.68\cdot10^{-4}$
	\end{tabular}
\end{center}
	\caption{Numerical results for full-order (in $y$) power series \eqref{eq:solpsa} corresponding 
	to Example 1.}
	\label{tab:numf}
\end{table}

Table~\ref{tab:numr} shows the corresponding results as Table~\ref{tab:numf} but with reduced 
order $N_y=2$ for $k(x,\xi,y)$ in \eqref{eq:kpsay}. No truncation has been applied to the 
parameter approximations \eqref{eq:ppca}, although the parameters other than $q$ are inherently 
low-order polynomials in $y$, so that the higher powers of $y$ do not appear in the Taylor series 
approximations anyway. In fact, from \eqref{eq:kc} and \eqref{eq:kcbca} if follows that, for the 
parameters of Example 1, only powers of $y$ of order $N_y=2$ appear in \eqref{eq:kpsay}. Thus, 
$N_y=2$ in the present case provides an exact approximation.\footnote{In this case, we, in fact, 
know a priori that $N_y=2$ is exact, because the closed-form 
	solution \eqref{eq:exkc} is a second-order polynomial in~$y$.}
Consequently, one can see that the number of unknowns in Table~\ref{tab:numr} grows 
much slower than in Table~\ref{tab:numf}, while the number of equations grows slower as well 
(although not 
as drastically). These lower numbers of unknowns and equations also reflect the slower increase 
in 
computational time as $N$ increases. When it comes to accuracy, both the residual of the 
least-squares fit and the maximal absolute error to the exact control gain decrease as $N$ 
increases, and notable improvements in accuracy are obtained for even $N$, similarly to 
Table~\ref{tab:numf}. Since the reduced-order power-series \eqref{eq:kpsay} is different from the 
full-order one, there are some discrepancies between the residuals and maximal errors when $N$ 
is small. Nevertheless, both approximations appear to converge to the exact solution at the same 
rate as $N$ increases.
\begin{table}[!htb]
	\begin{center}
		\begin{tabular}{c | c c c c c} 
			N & \#K & \#Eq & C. time & $\|\mathbf{Ax}-\mathbf{b}\|_2$ & max error \\ 
			\hline \\ [-11pt]
		12 & $326$ & $862$ & $7.23$ s. & $1.93$ & $16.17$ \\
		13 & $379$ & $999$ & $8.85$ s. & $0.673$ & $8.32$ \\
		14 & $436$ & $1146$ & $10.47$ s. & $0.210$ & $0.510$ \\
		15 & $497$ & $1303$ & $12.61$ s. & $5.55\cdot10^{-2}$  & $0.658$ \\
		16 & $562$ & $1470$ & $14.74$ s. & $1.34\cdot10^{-2}$  & $0.110$ \\
		17 & $631$ & $1647$ & $17.06$ s. & $3.10\cdot10^{-3}$  & $9.08\cdot10^{-2}$ \\
		18 & $704$ & $1834$ & $19.66$ s. & $6.84\cdot10^{-4}$  & $7.27\cdot10^{-3}$ \\
		19 & $781$ & $2031$ & $22.96$ s. & $1.42\cdot10^{-4}$  & $7.61\cdot10^{-3}$ \\
		20 & $862$ & $2238$ & $26.07$ s. & $2.82\cdot10^{-5}$  & $5.68\cdot10^{-4}$
		\end{tabular}
	\end{center}
	\caption{Numerical results of computation of reduced-order (in $y$) power series 
	\eqref{eq:kpsay} with $N_y=2$ corresponding to Example~1.}
	\label{tab:numr}
\end{table}

As the parameter $q$ only appears in an integral in the boundary condition \eqref{eq:kcbcb}, it 
may not in fact be necessary to approximate it by a power series, as long as the integral can be 
evaluated in closed form. Considering that $k$ in the integral is approximated by 
power series (and $\lambda$ is constant), the integral can generally be evaluated in closed form 
by using integration by parts, 
as long as the integral of $q$ over the unit interval can be computed in closed form. This is 
demonstrated in Table~\ref{tab:numrq}, which is analogous to Table~\ref{tab:numr}, except that 
the exact expression $q(y) = \cos(2\pi y)$ is used in the computations instead of a power series 
approximation. The data of Table~\ref{tab:numrq} is consistent with the data of 
Table~\ref{tab:numr}, except that for larger values of $N$ the maximal error to the exact solution 
is one order of magnitude smaller than with the power series approximation. Considering that this 
can be achieved in virtually the same computational time, it may be preferable to use the exact 
expression of $q$ in the computations (when possible) instead of a power series approximation.

\begin{table}[!htb]
	\begin{center}
		\begin{tabular}{c | c c c c c} 
			N & \#K & \#Eq & C. time & $\|\mathbf{Ax}-\mathbf{b}\|_2$ & max error \\ 
			\hline \\ [-11pt]
			12 & $326$ & $862$ & $5.27$ s. & $2.04$ & $13.14$ \\
			13 & $379$ & $999$ & $6.66$ s. & $0.712$ & $5.34$ \\
			14 & $436$ & $1146$ & $7.80$ s. & $0.207$ & $1.16$ \\
			15 & $497$ & $1303$ & $13.31$ s. & $5.44\cdot10^{-2}$  & $0.159$ \\
			16 & $562$ & $1470$ & $16.17$ s. & $1.34\cdot10^{-2}$  & $2.44\cdot 10^{-2}$ \\
			17 & $631$ & $1647$ & $18.94$ s. & $3.11\cdot10^{-3}$  & $4.16\cdot10^{-3}$ \\
			18 & $704$ & $1834$ & $19.37$ s. & $6.84\cdot10^{-4}$  & $7.42\cdot10^{-4}$ \\
			19 & $781$ & $2031$ & $23.23$ s. & $1.42\cdot10^{-4}$  & $1.32\cdot10^{-4}$ \\
			20 & $862$ & $2238$ & $27.80$ s. & $2.82\cdot10^{-5}$  & $2.27\cdot10^{-5}$
		\end{tabular}
	\end{center}
	\caption{Numerical results using the exact $q$ function, of reduced-order (in $y$) power series 
	\eqref{eq:kpsay} 
	with $N_y=2$ corresponding to Example~1.}
	\label{tab:numrq}
\end{table}

\section{Application to Stabilization of a Large-Scale System} \label{sec:ex2}

\subsection{Parameters and Their Continuum Approximation} \label{ex2:param}

\begin{figure}[htbp]
	\includegraphics[width=\columnwidth]{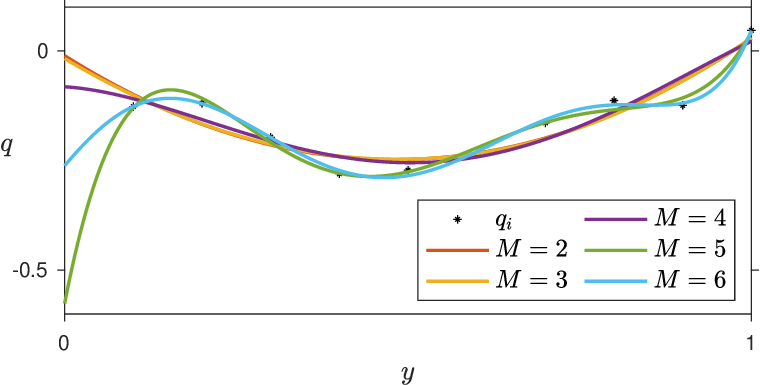}
	\caption{Polynomial fits of order $M=2,\ldots,6$ to the $q_i$ data.}
	\label{fig:pfex2}
\end{figure}

We consider a large-scale system of $n+1$ hyperbolic PDEs shown in \ref{app:n+1} 
(see \cite{DiMVaz13}) with parameters
\begin{subequations}
	\label{eq:ex2n+1}%
	\begin{align}
		\lambda_i(x) & = 1, \\
		\mu(x) & = 1, \\
		\sigma_{i,j}(x) & = x^3(x+1)\left(\frac{i}{n}-1\right)\left(\frac{j}{n}-1\right), \\	
		W_i(x) & = 2x(x+1)\frac{i}{n}, \\
		\theta_i(x) & = -70x\frac{i}{n}\left(\frac{i}{n}-1\right),
	\end{align}%
\end{subequations}
for $i,j = 1,2,\dots,n$, where we take $n=10$, and $q = 
[-0.127, -0.119, -0.197, -0.28, -0.272, -0.235, -0.164$, \\ $ -0.113,-0.124, 0.047]$\footnote{The 
values of $q$ are 
generated by sampling the polynomial $y(y-1)$ at $y=0.1,0.2,\ldots,1$ and then perturbing the 
data with evenly-distributed random numbers from the open interval $(-0.05,0.05)$.}. As 
discussed in Section 
\ref{sec:appr}, when $n$ is sufficiently large, such $n+1$ system can be approximated by a 
corresponding ensemble of linear PDEs, shown in \ref{app:n+1} (see \cite{AllKrs23arxiv, 
HumBek24arxiv}). This also implies that the solution 
to the large-scale kernel equations \eqref{eq:kn}, \eqref{eq:knbc} can be approximated by solving 
the corresponding continuum kernel equations \eqref{eq:kc}, \eqref{eq:kcbc} (see 
\cite{HumBek24arxiv} for details). In order to do this, 
we need to first construct the continuum parameters such that the relations \eqref{eq:afn1} are 
satisfied. Due to the structure of the parameters \eqref{eq:ex2n+1}, this can be achieved by 
replacing $i/n$ and $j/n$ by $y$ and $\eta$, respectively, which results in continuum parameters 
as
\begin{subequations}
	\label{eq:ex2}%
	\begin{align}
		\lambda(x,y) & = 1, \\
		\mu(x) & = 1, \\
		\sigma(x,y,\eta) & = x^3(x+1)(y-1)(\eta-1), \\	
		W(x, y) & = 2x(x+1)y, \\
		\theta(x,y) & = -70xy(y-1).
	\end{align}%
\end{subequations}
For the $q_i$ data, we construct a continuum approximation by utilizing the \verb+polyfit+ routine 
in MATLAB. We try polynomials of orders $M=2,\ldots,6$, which are illustrated in 
Figure~\ref{fig:pfex2}. Note that we have chosen to locate the $q_i$ data at points $y = i/n$ 
for $i=1,\ldots,n$, meaning that the leftmost data point is located at $1/n = 
1/10$.\footnote{One 
may just as well locate the points differently, e.g., at $(i-1)/n$ for $i=1,\ldots,n$ (see 
\cite[Rem. 4.4]{HumBek24arxiv} for details).} It can be seen 
that the fits of orders $M=2,3,4$ do not notably differ from one 
another, while the higher-order ones for $M=5,6$ provide a different, nevertheless mutually 
similar fit. We proceed with the lowest order fit $M=2$ to the power series computations, as that 
already seems to provide a sufficiently accurate approximation for the $q_i$ data. As we 
demonstrate at the end  of the next section (see Table \ref{tab:Mcomp}), this choice has only a 
marginal effect to the numerical results presented in the next section, as all the polynomial fits 
presented in Figure~\ref{fig:pfex2} are close to each other in the $L^2$ sense.

\subsection{Computation of Continuum Kernels with Power Series} \label{ex:num}

We note first that the continuum parameters constructed in Section~\ref{ex2:param} are 
separable, low-order polynomials in the spatial variables. However, as $\theta_y(y) = 
y(y-1)$ and $\sigma_y(y) = \sigma_{\eta}(y)  = y-1$, we have $\displaystyle 
\frac{\sigma_{\eta}(y)}{\theta_y(y)} = \frac{1}{y}$ 
and $ \displaystyle
\int\limits_0^1\sigma_y(\eta)\theta_y(\eta)d\eta = \frac{1}{12}$,
so that there is no constant $c_y$ satisfying \eqref{eq:cy}, and hence,  
Proposition~\ref{prop:esol} is not applicable for finding the solution to \eqref{eq:kc}, 
\eqref{eq:kcbc} in closed form. Instead, we employ the power series approach to solve the 
continuum kernel equations. Note that we were not able to find a closed-form solution to 
\eqref{eq:kn}, \eqref{eq:knbc} for parameters \eqref{eq:ex2n+1} either.

Considering that the parameters are low-order polynomials, one may guess that a power series 
approximation of the same order is sufficient for solving the equations. Hence, we initialize our 
computation of the kernel equations, with power series approximation of order $N=6$, as that is 
the highest 
order term appearing in the parameters shown in \eqref{eq:ex2}, namely $x^4y\eta$ in 
$\sigma(x,y,\eta)$. The results are shown in Table~\ref{tab:num2f}, which shows the data 
corresponding to Table~\ref{tab:numf} for this example. However, as we do not know the exact 
solution to the 
kernel equations \eqref{eq:kc}, \eqref{eq:kcbc} for parameters \eqref{eq:ex2}, instead of the 
maximal error to the exact solution, the penultimate column shows the 
maximal difference $d_{n+1}$ to the solution of the large-scale kernel equations, which we have 
computed 
based on a finite-difference approximation of \eqref{eq:kn}, \eqref{eq:knbc}. Moreover, the last 
column shows the maximal difference $d_{N-1}$ to the solution obtained with power series of 
order $N-1$. 
The last two columns imply that the power series solutions to \eqref{eq:kc}, \eqref{eq:kcbc} 
converge as $N$ increases. However, as $n=10$, there persists an approximation error between 
the 
solutions to \eqref{eq:kn}, \eqref{eq:knbc} and \eqref{eq:kc}, \eqref{eq:kcbc}, which would tend to 
zero as $n\to\infty$ (see \cite{HumBek24arxiv}). Moreover, the second-order polynomial fit to the 
$q_i$ data also induces some (additional) approximation errors in this case.

\begin{figure}[htbp]
	\includegraphics[width=\columnwidth]{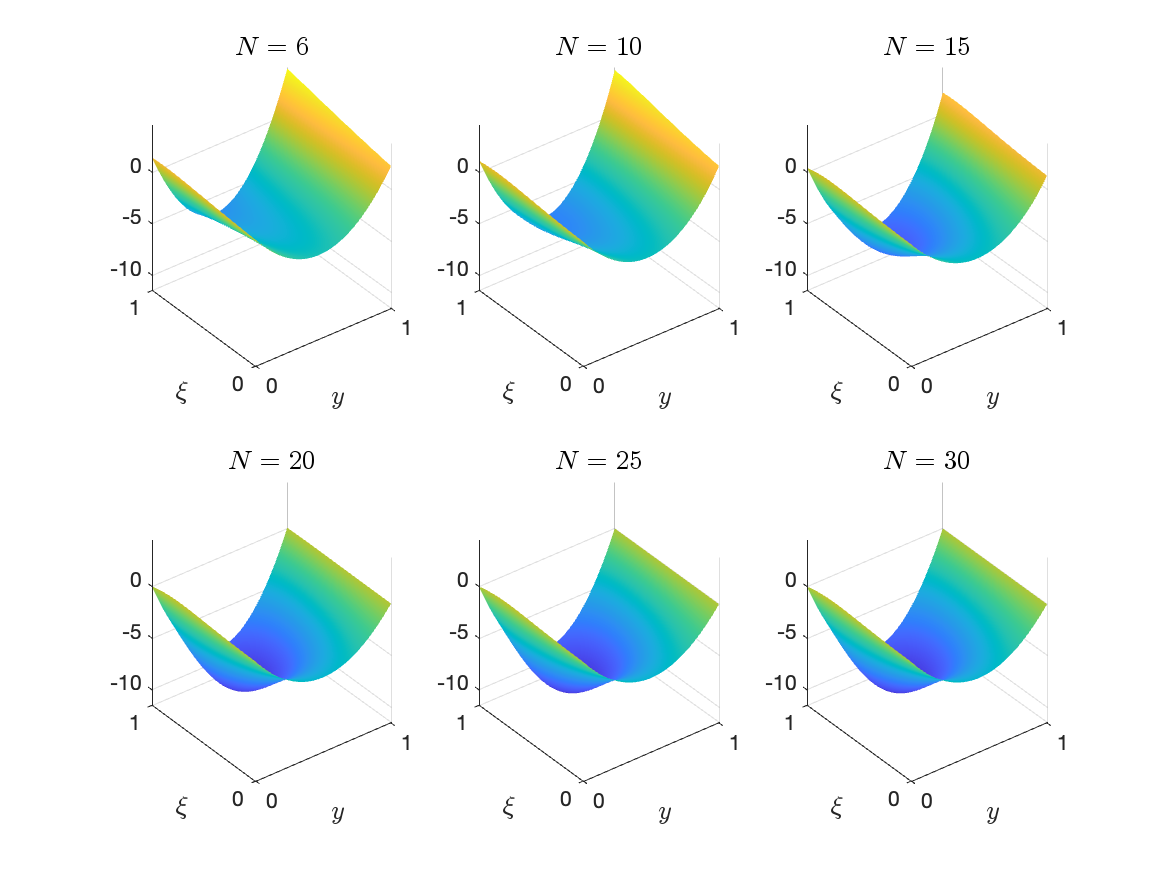}
	\caption{The control gain $k(1,\xi,y)$ approximated by full-order (in $y$) power series 
	\eqref{eq:kpsa} for $N = 
	6,10,15,20,25,30$.}
	\label{fig:ex2}
\end{figure}

The obtained 
approximate solutions \eqref{eq:kpsa} evaluated at $x=1$ are displayed in Figure~\ref{fig:ex2}, 
where visually it is 
difficult to distinguish the solutions after $N=20$, further supporting the conclusion of 
convergent power series approximations. However, it is interesting to note that a much higher 
order approximation is needed for the solution than the powers of the parameters \eqref{eq:ex2}, 
which, in this 
case, are represented exactly by a Taylor series of order six. In fact, the convergence rate of the 
power series approximation appears to be slower than in Example 1 (based on the values of the 
least squares fit  $\|\mathbf{Ax}-\mathbf{b}\|_2$), where the approximation 
accuracy is good already at the order of around $N=15$, even though the parameters are not 
polynomials. Regardless, in Example 1 both the parameters and the solution \eqref{eq:exkc} can 
be approximated at a satisfactory accuracy with this order of approximation. In general, 
as we may have no a priori information about the form of the solution, the convergence rate 
and approximation accuracy need to be assessed based on the error
of the least squares fit $\|\mathbf{Ax}-\mathbf{b}\|_2$ as $N$ increases.\footnote{Although here 
we derive such a convergence rate proxy based on numerical experiments, in principle, an 
upper bound for the convergence rate of the power series could be derived analytically. 
However, this would require estimating the derivatives of $(k,\bar{k})$ up to arbitrary order, 
which, based on the discussion at the end of Section~\ref{pow:alg}, may be computationally 
laborious.  Moreover, it is, in general, expected that such upper bounds may be conservative for 
practical computations.}

\begin{table}[!htb]
	\begin{center}
		\resizebox{.99\columnwidth}{!}{
		\begin{tabular}{c | c c c c c c} 
			N & \#K & \#Eq & C. time & $\|\mathbf{Ax}-\mathbf{b}\|_2$ & $d_{n+1}$ & $d_{N-1}$ \\ 
			\hline \\ [-11pt]
			6 & $112$ & $221$ & $1.50$ s. & $4.64$ & $5.90$ & $-$ \\
			10 & $352$ & $533$ & $4.74$ s. & $3.46$ & $4.78$ & $0.161$ \\
			15 & $952$ & $1223$ & $15.09$ s. & $2.07$ & $2.86$ & $0.127$\\
			20 & $2002$ & $2363$ & $36.70$ s. & $0.414$  & $1.15$ & $0.184$ \\
			25 & $3627$ & $4078$ & $89.65$ s. & $2.6\cdot10^{-2}$  & $1.09$ & $2.1\cdot 
			10^{-2}$ \\
			30 & $5952$ & $6493$ & $212.29$ s. & $9.3\cdot10^{-4}$  & $1.09$ & $1.6\cdot 
			10^{-5}$ 
		\end{tabular}
	}
	\end{center}
	\caption{Numerical results for full-order (in $y$) power series approximation \eqref{eq:solpsa} 
	corresponding to parameters \eqref{eq:ex2}.}
	\label{tab:num2f}
\end{table}

\begin{figure}[htbp]
	\includegraphics[width=\columnwidth]{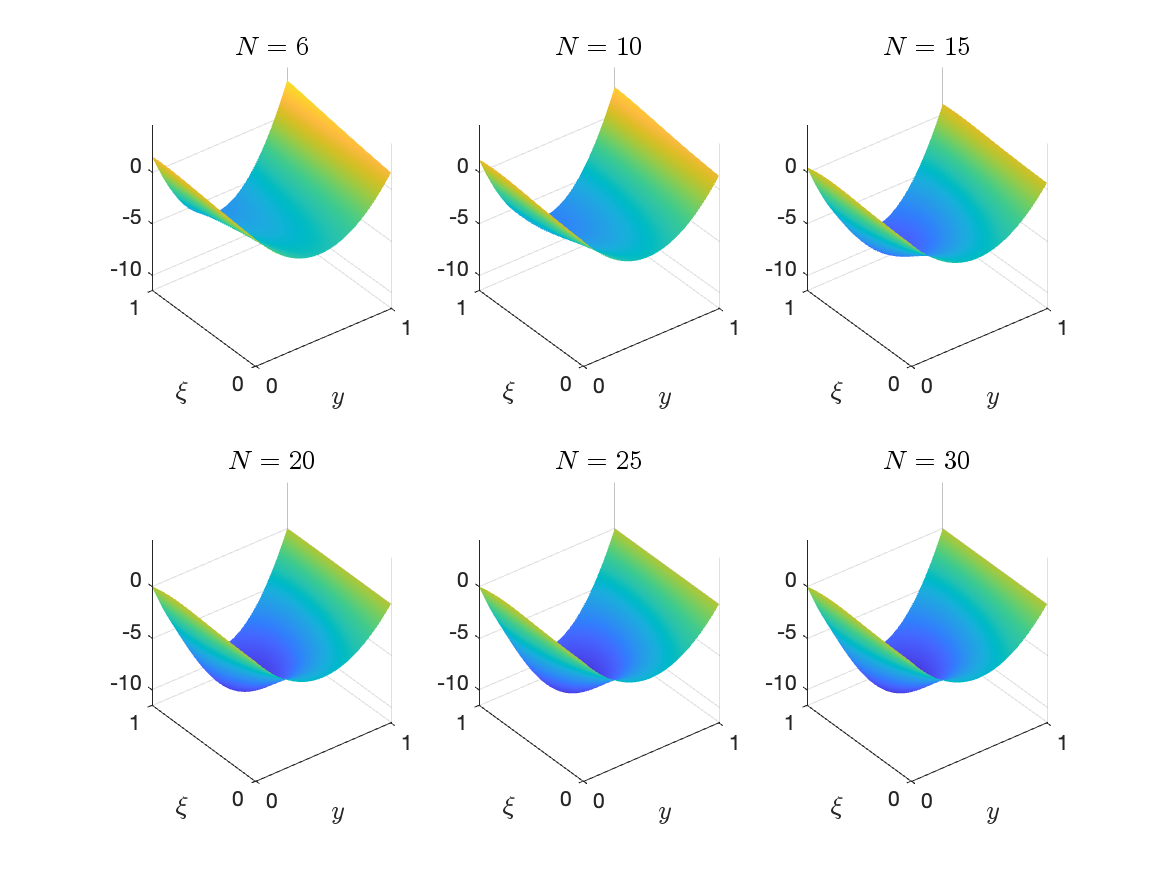}
	\caption{The control gain $k(1,\xi,y)$ approximated by reduced-order (in $y$) power series 
	\eqref{eq:kpsay}  for 
	$N = 6,10,15,20,25,30$ and $N_y=2$.}
	\label{fig:ex2r}
\end{figure}

We employ also here the reduced-order approximation in $y$. Considering that the highest power 
of $y$ appearing in the parameters is two, we repeat the computations of 
Table~\ref{tab:num2f} with reduced order $N_y = 2$ in $y$. Similarly to what we saw in 
Example~1, 
the results shown in Table~\ref{tab:num2r} for the reduced-order (in $y$) power series are 
analogous to the full-order results in terms of accuracy, apart from some minor discrepancies with 
smaller values of $N$. The obtained control gains $k(1,\xi,y)$ for different $N$ are displayed in 
Figure~\ref{fig:ex2r}, 
where one can notice some differences for smaller values of $N$ as compared to 
Figure~\ref{fig:ex2}. For larger values of $N$ (bottom plot of Figure~\ref{fig:ex2}), the two figures 
are 
virtually identical. 
Thus, the reduced-order power series provides as good of an approximation (when $N$ is 
sufficiently large) as the full-order power series, with the notable benefit that the 
computational complexity (number of unknowns/equations and computational time) 
increases at a slower rate with respect to $N$ than in the full-order power series approximation, 
as already reported in Table~\ref{tab:ccc}.

\begin{table}[!htb]
	\begin{center}
		\resizebox{.99\columnwidth}{!}{
		\begin{tabular}{c | c c c c c c} 
			N & \#K & \#Eq & C. time & $\|\mathbf{Ax}-\mathbf{b}\|_2$ & $d_{n+1}$ & $d_{N-1}$ \\ 
			\hline \\ [-11pt]
			6 & $92$ & $201$ & $1.37$ s. & $4.78$ & $5.95$ & $-$ \\
			10 & $232$ & $413$ & $3.21$ s. & $3.62$ & $4.83$ & $0.173$ \\
			15 & $497$ & $768$ & $7.37$ s. & $2.14$ & $2.82$ & $0.138$ \\
			20 & $862$ & $1223$ & $13.03$ s. & $0.417$  & $1.15$ & $0.180$ \\
			25 & $1327$ & $1778$ & $21.52$ s. & $2.7\cdot10^{-2}$  & $1.09$ & $2.0\cdot 
			10^{-2}$ \\
			30 & $1892$ & $2433$ & $33.98$ s. & $9.3\cdot10^{-4}$  & $1.09$ &  $1.9\cdot 
			10^{-5}$ 
		\end{tabular}
	}
	\end{center}
	\caption{Numerical results for reduced-order (in $y$) power series approximation 
	\eqref{eq:kpsay} with $N_y=2$ corresponding to the parameters \eqref{eq:ex2}.}
	\label{tab:num2r}
\end{table}

To conclude this section, we justify our choice of $M=2$ for the order of the polynomial fit to the 
$q_i$ data. The results are displayed in Table \ref{tab:Mcomp} for full-order power series (in $y$) 
of order $N=20$, which show that the different orders of polynomial fit to the $q_i$ data have a 
marginal effect on the results. Thus, choosing the lowest order fit of $M=2$ is reasonable, as that 
results in the simplest continuum approximation.
\begin{table}[!htb]
	\begin{center}
			\begin{tabular}{c | c c c c c}
				M & \#K & \#Eq & C. time & $\|\mathbf{Ax}-\mathbf{b}\|_2$ & $d_{n+1}$ \\ 
				\hline \\ [-11pt]
				2 & $2002$ & $2363$ & $34.96$ s. & $0.4136$ & $1.1484$ \\
				3 & $2002$ & $2363$ & $35.49$ s. & $0.4136$ & $1.1484$ \\
				4 & $2002$ & $2363$ & $35.78$ s. & $0.4136$ & $1.1486$ \\
				5 & $2002$ & $2363$ & $36.13$ s. & $0.4140$  & $1.1497$ \\
				6 & $2002$ & $2363$ & $36.44$ s. & $0.4138$  & $1.1490$ 
			\end{tabular}
	\end{center}
	\caption{Numerical results for full-order (in $y$) power series approximation 
	\eqref{eq:solpsa} of order $N=20$ corresponding to parameters \eqref{eq:ex2} with 
	polynomials of order $M = 2,\ldots,6$ fitted to the $q_i$ as displayed in Figure~\ref{fig:pfex2}.	
	\label{tab:Mcomp}}
\end{table}

\subsection{Stabilization Using Power Series-Based Approximate Continuum Kernels}

\begin{figure}[htbp]
	\includegraphics[width=\columnwidth]{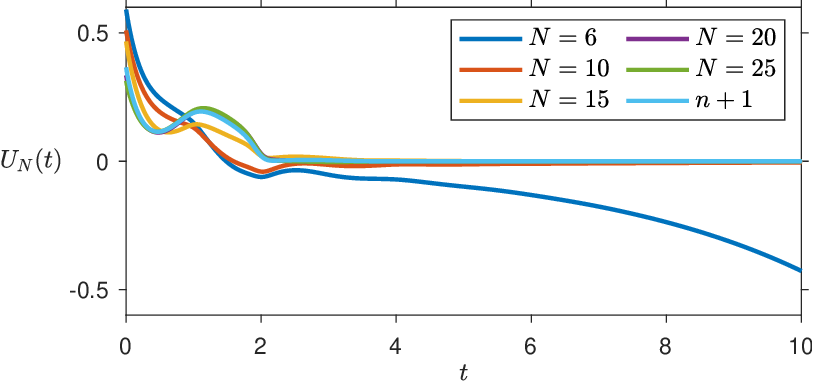}
	\caption{The controls obtained with approximate continuum kernels computed with full-order 
	(in $y$) power series \eqref{eq:solpsa} of order $N = 6,10,15,20,25$ and by solving 
	\eqref{eq:kn}, 	\eqref{eq:knbc} for the $n+1$ kernels.}
	\label{fig:ex2U}
\end{figure}

\begin{figure}[htbp]
	\includegraphics[width=\columnwidth]{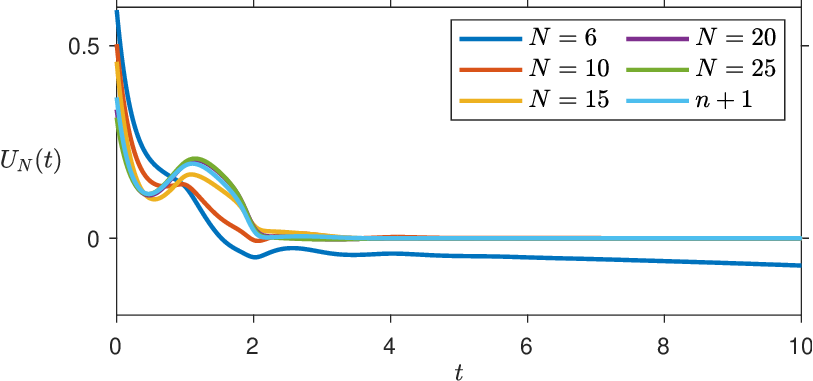}
	\caption{The controls obtained with approximate continuum kernels computed with 
	reduced-order (in $y$) power series \eqref{eq:kpsay} of order $N = 6,10,15,20,25$ and 
	$N_y=2$ and by solving \eqref{eq:kn}, \eqref{eq:knbc} for the $n+1$ kernels.}
	\label{fig:ex2Ur}
\end{figure}

We simulate the $n+1$ system with parameters \eqref{eq:ex2n+1}. The $n+1$ hyperbolic PDE 
system, is 
approximated  using finite differences with $256$ grid 
points in $x$. For the backstepping control law, we employ the continuum kernels computed in 
Section~\ref{ex:num} for different orders of approximation $N$. Note that with parameters 
\eqref{eq:ex2}, it is verified in simulation that the open-loop system is unstable. The 
employed hardware and software are the same as in Section 5, and the data displayed in Figures 
\ref{fig:ex2U} and \ref{fig:ex2Ur} takes about 14 seconds to compute by solving the 
finite-difference approximation of the closed-loop system with the \texttt{ode45} solver in 
MATLAB.

Figure~\ref{fig:ex2U} shows the control inputs corresponding to the continuum kernels computed 
using the full-order (in $y$) power series of order $N = 6,10,15,20,25$ and the 
solution\footnote{The 
solution is computed based on a finite-difference approximation of \eqref{eq:kn}, 
\eqref{eq:knbc}.} to the $n+1$ kernel 
equations \eqref{eq:kn}, \eqref{eq:knbc}. The controls start to diverge for 
$N=6$, implying that the controls fail to stabilize the closed-loop system. This is expected 
based on the numerical experiments of Section~\ref{ex:num}, as the accuracy obtained with such 
a low-order approximation is not adequate. As the approximation order $N$ increases, 
one can see that the controls tend to zero, implying that, for larger values of $N$, the 
controls are stabilizing. For the lower orders $N=10$ and $N=15$ the controls are distinguishable, 
and the decay rate is relatively slow, whereas for $N=20$ and $N =25$ the controls are virtually
identical and converge to zero after two seconds, implying that all the state 
components have converged to zero by that time as well. This is consistent with the control gains 
displayed in Figure~\ref{fig:ex2}, which are distinctly different for $N=5,10,15$ and virtually 
identical 
for $N=20,25,30$. Finally, we see that the approximate controls for $N=20$ and $N=25$ 
are also very close to the controls computed based on the solution to the $n+1$ kernel equations 
\eqref{eq:kn},~\eqref{eq:knbc}.

Figure~\ref{fig:ex2Ur} shows the corresponding results as Figure~\ref{fig:ex2U} but for the 
reduced-order (in $y$) power series approximation \eqref{eq:kpsay} for the continuum kernel. 
The same 
comments apply to Figure~\ref{fig:ex2Ur} as to Figure~\ref{fig:ex2U}. A slight difference can be 
observed for the low-order approximations, where the case $N=6$ diverges at a slower rate than 
in Figure~\ref{fig:ex2U} (the closed-loop system is still unstable), while the case $N=10$ appears 
to 
converge to zero at a faster rate. Regardless, for larger values of $N$, the controls are 
virtually indistinguishable from one another, as well as from the ones displayed in 
Figure~\ref{fig:ex2U}.
This is consistent with the control gains in Figure~\ref{fig:ex2r}. We note here that as \cite[Thm 
4.1]{HumBek24arxiv} predicts, for sufficiently large $n$ ($n=10$ in the present case), the 
corresponding feedback law obtained employing the continuum control kernel $k(1,\xi,y)$ is 
stabilizing, provided that $k(1,\xi,y)$ is approximated at a sufficient accuracy via the power series. 
In turn, the latter is established in Theorem~\ref{thm:pswp}. In conclusion, stabilization with 
performance close to the one corresponding to exact control kernels is achieved, with 
computational complexity of order $\mathcal{O}(N^2)$ that does not grow with $n$, in contrast 
to the computation of the exact control kernels (see Table~\ref{tab:ccc}). We note that, the 
computational complexity of solving \eqref{eq:kn}, \eqref{eq:knbc} using a finite-difference 
approximation grows (roughly linearly) with $n$, as shown in \cite[Figure 6]{HumBek24arxiv}.

\section{Conclusions and Future Work} \label{sec:conc}

In this paper, we provided computational tools for construction of backstepping-based stabilizing 
kernels for continua of hyperbolic PDEs, 
and thus, also for large-scale PDE systems. These tools include both explicit kernels construction 
(even though the class of systems for which this is possible 
may be restrictive) and approximation via power series. We demonstrated the accuracy and 
efficiency of these tools on 
two numerical examples. In the first, we derived the closed-form solution to 
continuum kernel equations and demonstrated the convergence of the power series 
approximation to the closed-form solution as $N$ increases. In the second numerical example, we 
could not derive the solution in closed form, and thus, we employed only the power series 
approach (and reduced-order power series approach) to 
solve the continuum kernel equations. We further demonstrated that the computed continuum 
kernels provided a stabilizing feedback law for the corresponding large-scale system of linear 
PDEs, with performance comparable to the case of employing the exact, large-scale control 
kernels.

For future work, the 
efficacy of Algorithm~\ref{alg:ps} could, potentially, be substantially improved by adopting the 
double vector-matrix framework from \cite[Sect. II.B]{LinVaz24arxiv}. However, this has to be 
modified to fit the triple power series employed here, potentially requiring a triple 
vector-matrix framework, which may not be straightforward to implement. Another topic for future 
research is to derive analytically the convergence rate of the power 
series representation for the solution to the continuum kernel equations, thus providing explicit 
(although potentially conservative) estimates of the lowest power required for an accurate, power 
series-based approximation (given a desired approximation error). However, as discussed at 
the end of Section~\ref{pow:alg}, derivation of such bounds is of the same computational 
complexity $\mathcal{O}(N^3)$ as the proposed power series procedure for actually computing 
the kernels. Hence, the 
least-squares error of Algorithm~ \ref{alg:ps} may be the most practical option for assessing the 
accuracy of the computed solution.

\appendix

\section*{Appendix}

\section{Proof of Proposition~\ref{prop:esol}} \label{app:esol}

The solution stated in Proposition~\ref{prop:esol} can be constructed by looking for a separable 
solution to \eqref{eq:kc}, \eqref{eq:kcbc}, 
i.e., $k(x,\xi,y)  = k_x(x)k_\xi(\xi)k_y(y)$ and $\bar{k}(x,\xi)= \bar{k}_x(x)\bar{k}_\xi(\xi)$. After 
dividing the equations \eqref{eq:kc} by the solution, we get
\begin{subequations}
	\label{eq:cker}%
	\begin{align}
		\mu \frac{k_x'(x)}{k_x(x)} -
		\lambda(y)\frac{k_{\xi}'(\xi)}{k_{\xi}(\xi)} -
		\frac{\theta_x(\xi)\theta_y(y)\bar{k}_x(x)\bar{k}_{\xi}(\xi)}{k_x(x)k_y(y)k_{\xi}(\xi)}
		& = \nonumber \\
		\frac{\sigma_x(\xi)\sigma_\eta(y)}{k_y(y)}\int\limits_0^1\sigma_{y}(\eta)k_y(\eta)d\eta, & 
		\label{eq:ckera}
		\\
		\mu\frac{\bar{k}'_x(x)}{\bar{k}(x)} + \mu
		\frac{\bar{k}'_{\xi}(\xi)}{\bar{k}_{\xi}(\xi)}
		= 
		\frac{W_x(\xi)k_x(x)k_{\xi}(\xi)}{\bar{k}_x(x)\bar{k}_{\xi}(\xi)}\int\limits_0^1W_y(y)k_y(y)dy, & 
		\label{eq:ckerb}
	\end{align}
\end{subequations}
with boundary conditions 
\begin{subequations}
	\label{eq:ckerbc}%
	\begin{align}
		k_x(x)k_{\xi}(x)k_y(y)
		& = -\frac{\theta_x(x)\theta_y(y)}{\lambda(y) + \mu}, \label{eq:ckerbca} \\
		\mu\bar{k}_x(x)\bar{k}_{\xi}(0) \
		& = k_x(x)k_{\xi}(0)\int\limits_0^1\lambda(y)q(y)k_y(y)dy. \label{eq:ckerbcb}
	\end{align}
\end{subequations}
As the boundary conditions have to hold for all $x,y \in [0,1]$, the second boundary condition 
implies that $\bar{k}_x(x)/k_x(x)$ has to be constant. Inserting this into \eqref{eq:ckera}, we also 
get that $\mu k_x'(x)/k_x(x)$ has to be constant, meaning that we can set 
$k_x(x) = \exp \left( \frac{c_x}{\mu}x\right) = \bar{k}_x(x)$ for some $c_x$, which 
eventually 
turns 
out to be given by \eqref{eq:cx}.

The first boundary condition \eqref{eq:ckerbca} gives 
\begin{equation} \label{eq:tmp1}
	(\lambda(y)+\mu)\frac{k_y(y)}{\theta_y(y)} = -
	\frac{\theta_x(x)}{k_x(x)k_{\xi}(x)},
\end{equation}
meaning that we can set 
\begin{equation}
	k_y(y) = c_1\frac{\theta_y(y)}{\lambda(y) + \mu},
\end{equation}
for some $c_1 \neq 0$, and further get 
\begin{equation} \label{eq:kxi}
k_{\xi}(x) = -\frac{1}{c_1}\frac{\theta_x(x)}{k_x(x)} =
	-\frac{1}{c_1}\theta_x(x)\exp \left( -\frac{c_x}{\mu}x \right),
\end{equation}
which, in combination with \eqref{eq:tmp1}, gives \eqref{eq:ksepa}. Now, substituting 
\eqref{eq:ksepa} and $\bar{k}_x(x) = \exp \left( \frac{c_x}{\mu}x\right)$ to \eqref{eq:ckera}, 
we 
obtain
\begin{align}
	c_x - \lambda(y) \left(\frac{\theta_x'(\xi)}{\theta_x(\xi)} -
	\frac{c_x}{\mu}\right) + \left( \lambda(y)+\mu \right)\exp \left(
	\frac{c_x}{\mu}\xi \right)\bar{k}_{\xi}(\xi) & = \nonumber \\ 
	\sigma_x(\xi)(\lambda(y)+\mu)\frac{\sigma_\eta(y)}{\theta_y(y)}\int\limits_0^1 
	\frac{\sigma_y(\eta)\theta_y(\eta)}{\lambda(\eta) + \mu}d\eta,
\end{align}
which gives $\bar{k}_{\xi}(\xi) = f(\xi)\exp \left( -\frac{c_x}{\mu}\xi \right)$ with $f$ given in 
\eqref{eq:f}. Now, multiplying \eqref{eq:ckerb} by
$\bar{k}_{\xi}(\xi)$ and using 
\begin{equation}
\bar{k}_\xi'(\xi) = -\frac{c_x}{\mu}\bar{k}_\xi(\xi) + f'(\xi) \exp \left( -\frac{c_x}{\mu}\xi \right),
\end{equation}
we get
\begin{align}
	f'(\xi) \exp \left( -\frac{c_x}{\mu}\xi \right) & = \nonumber \\
 W_x(\xi)\theta_x(\xi)\exp \left( -\frac{c_x}{\mu}\xi 
	\right) \int\limits_0^1\frac{W_y(y)\theta_y(y)}{\lambda(y) + \mu}dy, &
\end{align}
which is satisfied as \eqref{eq:fd} holds by assumption.
Finally, from the second boundary condition \eqref{eq:ckerbcb} we get
\begin{equation}
	f(0) = 
	-\frac{\theta_x(0)}{\mu}\int\limits_0^1\lambda(y)q(y)\frac{\theta_y(y)}{\lambda(y) + \mu}dy,
\end{equation}
which, when combined with \eqref{eq:f}, gives the stated expression \eqref{eq:cx} for $c_x$. This 
concludes the proof.

\section{Large-Scale System of $n+1$ Hyperbolic PDEs and Its Continuum Approximation} 
\label{app:n+1}

By an $n+1$ system, or a large-scale system, we refer to the following system of $n+1$ linear 
hyperbolic PDEs
\begin{subequations}
	\label{eq:n+1}%
	\begin{align}
		u_t^i(t,x) + \lambda_i(x)u_x^i(t,x)
		& = \frac{1}{n}\sum_{j=1}^n\sigma_{i,j}(x)u^j(t,x) \nonumber \\
		 & \qquad + W_i(x)v(t,x), \\
		v_t(t,x) - \mu(x)v_x(t,x)
		& =  \frac{1}{n}\sum_{j=1}^n \theta_j(x)u^j(t,x),
	\end{align}
\end{subequations}
with boundary conditions
\begin{equation}
	\label{eq:nbcuy}
	u^i(t,0)  = q_iv(t,0), \qquad v(t,1) = U(t), 
\end{equation}
for $i=1,2,\ldots,n$. It follows from \cite[Thm 3.2]{DiMVaz13} (see also \cite[Sect. 
II]{HumBek24arxiv}) that the system  \eqref{eq:n+1}, 
\eqref{eq:nbcuy} is exponentially stabilizable by a state-feedback law of the form 
\begin{equation}
	\label{eq:Un}
	U(t) = \int\limits_0^1 \left[  \frac{1}{n}\sum_{i=1}^nk^i(1,\xi)u^i(t,\xi) +
	k^{n+1}(1,\xi)v(t,\xi) \right] d\xi,
\end{equation}
where $\left(k^i\right)_{i=1}^{n+1}$ is the solution to \eqref{eq:kn}, \eqref{eq:knbc}. As already 
discussed in Section~\ref{sec:appr}, the $n+1$ kernel equations \eqref{eq:kn}, \eqref{eq:knbc} 
can be approximated by the corresponding continuum kernel equations 
\eqref{eq:kc}, \eqref{eq:kcbc}, provided that the continuum parameters are constructed 
appropriately (see \cite[Thm 4.1]{HumBek24arxiv}). For example, they can be constructed such 
that the relations in \eqref{eq:afn1} are satisfied. Under an appropriate 
continuum approximation of the parameters of the $n+1$ system, provided that $n$ is sufficiently 
large, the control gains in 
\eqref{eq:Un} can be replaced by the corresponding continuum gains given in  
\eqref{eq:kappr}, with preservation of exponential stability \cite[Thm 4.1]{HumBek24arxiv}.

The continuum kernel equations \eqref{eq:kc}, \eqref{eq:kcbc} first appeared within the 
framework of backstepping control design for an ensemble (or continuum) of hyperbolic PDEs of 
the form
\begin{subequations}
	\label{eq:inf}%
	\begin{align}
		u_t(t,x,y) + \lambda(x,y)u_x(t,x,y) \nonumber  & = \int\limits_0^1 
		\sigma(x,y,\eta)u(t,x,\eta)d\eta 
		\nonumber \\
		& \quad  + W(x,y)v(t,x),   \\
		v_t(t,x) - \mu(x)v_x(t,x)  & = \int\limits_0^1\theta(x,y)u(t,x,y)dy,
	\end{align}
\end{subequations}
with boundary conditions 
\begin{equation}
	\label{eq:cbcuy}
	u(t,0,y) = q(y)v(t,0), \qquad v(t,1) = U(t),
\end{equation}
for almost every $y\in[0,1]$, where $y$ is the ensemble variable. The stabilizing backstepping 
control law for the continuum system \eqref{eq:inf}, \eqref{eq:cbcuy} is given by \cite[Thm 
1]{AllKrs23arxiv}
\begin{equation}
	\label{eq:Uc}
	U(t) = \int\limits_0^1 \left[\int\limits_0^1 k(1,\xi,y)u(t,\xi,y)dy
	+ \bar{k}(1,\xi)v(t,\xi) \right]d\xi,
\end{equation}
where $(k,\bar{k})$ is the solution to \eqref{eq:kc}, \eqref{eq:kcbc}. In fact, for large $n$, the 
solutions to \eqref{eq:n+1}, \eqref{eq:nbcuy} converge to the solutions to \eqref{eq:inf}, 
\eqref{eq:cbcuy} as well, under appropriate conditions \cite[Thm 6.1]{HumBek24arxiv}.

\end{document}